\documentclass[onefignum,onetabnum]{siamart220329}

\usepackage{subfigure}
\usepackage{booktabs}
\usepackage{amsmath,amsxtra,amsfonts,amscd,amssymb,bm}
\usepackage{exscale}
\usepackage{geometry}
\usepackage{dcolumn}
\usepackage[shortlabels]{enumitem}
\usepackage{multirow}
\usepackage{psfrag}
\psfrag{SGS-FP}{\scalebox{.6}{\!\!\!SGS-FP}}
\psfrag{SGS-Syn}{\scalebox{.6}{\!\!\!SGS-PFP}}
\psfrag{SGS-MM}{\scalebox{.6}{\!\!\!SGS-PFP}}

\newtheorem{remark}{Remark}
\newcommand\bs[1]{\boldsymbol{#1}}

\makeatletter
\def\namedlabel#1#2{\begingroup
    #2%
    \def\@currentlabel{#2}%
    \phantomsection\label{#1}\endgroup
}
\makeatother

\begin{document}
\title{Symmetric Gauss-Seidel Method with a Preconditioned Fixed-Point Iteration for the Steady-State Boltzmann equation\thanks{Submitted to the editors.
\funding{The work of the first author was supported by the Academic Research Fund of the Ministry of Education of Singapore grant A-8000695-00-00. The work of the last author is partially supported under the NSF grant DMS-2409858, AFOSR grant FA9550-21-1-0358, and DOE grant DE-SC0023164.}}}

\author{
Zhenning Cai\thanks{Department of Mathematics, National University of Singapore, Singapore 119076 (\email{matcz@nus.edu.sg}).}
\and
Xiaoyu Dong\thanks{Engineering Systems and Design, Singapore University of Technology and Design, Singapore 487372 (\email{dongxiaoyu@lsec.cc.ac.cn}).}
\and
Jingwei Hu\thanks{Department of Applied Mathematics, University of Washington, Seattle, WA 98195, USA (\email{hujw@uw.edu}).}
}

\maketitle

\begin{abstract}
We introduce a numerical solver for the steady-state Boltzmann equation based on the symmetric Gauss-Seidel (SGS) method. To solve the nonlinear system on each grid cell derived from the SGS method, a fixed-point iteration preconditioned with its asymptotic limit is developed. The preconditioner only requires solving an algebraic system which is easy to implement and can speed up the convergence significantly especially in the case of small Knudsen numbers. Additionally, we couple our numerical scheme with the multigrid method to accelerate convergence. A variety of numerical experiments are carried out to illustrate the effectiveness of these methods.
\end{abstract}


\section{Introduction} 

Rarefied gases, which refer to gases at low densities or pressures, are crucial in various scientific and engineering fields, including aerospace engineering, vacuum technology, and plasma physics. It also plays an important role in the development of micro-electromechanical systems and nanotechnology, where gas behavior at small scales is of paramount importance. Knudsen number is one of the key parameters to characterize rarefied gases, which is the ratio of the molecular mean free path to a characteristic length scale. When the Knudsen number is large, the gas is considered rarefied, and molecular collisions and interactions become less significant. Compared to classical continuum fluid dynamics equations like the Navier-Stokes equations, Boltzmann equation provides a more accurate representation of rarefied gases, which is a mesoscopic kinetic model allowing us to bridge the gap between the macroscopic world of fluid dynamics and the microscopic world of individual gas molecules. However, compared with fluid equations, the Boltzmann equation contains an additional velocity variable, making the numerical simulation much more challenging. In this work, we focus on the numerical solver of the steady-state Boltzmann equation, where the fluid has reached a constant state that does not change with time.

In the gas kinetic theory, the distribution function $f(\bs{x}, \bs{v})$ is adopted to describe the gas states, with $\bs{x}$ and $\bs{v}$ being the spatial and velocity variables, respectively. The steady-state Boltzmann equation then reads
\begin{equation} \label{eq:Boltz}
\bs{v} \cdot \nabla_{\bs{x}} f(\bs{x},\bs{v}) = \frac{1}{\epsilon} \mathcal{Q}[f,f](\bs{x}, \bs{v}), \qquad \bs{x} \in \Omega \subset \mathbb{R}^d, \quad \bs{v} \in \mathbb{R}^d,
\end{equation}
where $\epsilon$ stands for Knudsen number and the operator $\mathcal{Q}$ describes the binary collision between gas molecules:
\begin{equation} \label{eq:collision}
\mathcal{Q}[f,f](\bs{x}, \bs{v}) = \int_{\mathbb{R}^d} \int_{\mathbb{S}^{d-1}} B(\bs{v}-\bs{v}_*,\bs{\omega}) \left[f(\bs{x}, \bs{v}'_*) f(\bs{x}, \bs{v}') - f(\bs{x}, \bs{v}_*) f(\bs{x}, \bs{v}) \right] \mathrm{d} \bs{\omega} \,\mathrm{d} \bs{v}_*.
\end{equation}
In the equation above, $\bs{v}$ and $\bs{v}_*$ denote the pre-collision velocities and the 
post-collision ones $\bs{v}'$ and $\bs{v}'_*$ are derived by
\begin{equation*}
\bs{v}' = \frac{\bs{v} + \bs{v}_*}{2} + \frac{|\bs{v} - \bs{v}_*|}{2} \bs{\omega}, \qquad
\bs{v}'_* = \frac{\bs{v} + \bs{v}_*}{2} - \frac{|\bs{v} - \bs{v}_*|}{2} \bs{\omega}.
\end{equation*}
The collision kernel $B(\cdot,\cdot)$ depends on the relative velocity $\bs{v} - \bs{v}_*$ and the impact angle $\bs{\omega}$, which is usually determined by the species of the gas.
The difficulty in solving the Boltzmann equation is two-fold. First, the distribution function includes six variables, and the collision term involves a five-dimensional integral, causing a high computational complexity in deterministic numerical methods. Second, an iterative method for this nonlinear equation is needed to achieve fast convergence. The complexity of the collision term can be alleviated by modeling the collision between particles using the Bhatnagar-Gross-Krook (BGK) approximation \cite{Bhatnagar1954}, where the collision term $\mathcal{Q}[f,f]$ is replaced by
\begin{equation} \label{eq:BGK_collis}
\mathcal{Q}^{\mathrm{BGK}}[f] = \nu \left( \mathcal{M}[f] - f \right),
\end{equation}
where $\nu$ is the collision frequency determining how quickly the distribution function approaches equilibrium, and $\mathcal{M}[f]$ represents the local Maxwellian defined by
\begin{equation} \label{eq:maxwellian}
\mathcal{M}[f] (\bs{x}, \bs{v}) = \frac{\rho(\bs{x})}{(\sqrt{2 \pi T(\bs{x})})^d} \exp \left( - \frac{|\bs{v} - \bs{U}(\bs{x})|^2}{2 T(\bs{x})} \right)
\end{equation}
in which $\rho, \bs{U}$ and $T$ are moments of the distribution function representing the density, velocity and thermal temperature of the gas, respectively. These moments are defined by
\begin{equation} \label{eq:moment}
\rho(\bs{x}) = \int_{\mathbb{R}^d} f(\bs{x},\bs{v}) \,\mathrm{d} \bs{v}, \quad
\bs{U}(\bs{x}) = \frac{1}{\rho(\bs{x})} \int_{\mathbb{R}^d} \bs{v} f(\bs{x},\bs{v}) \,\mathrm{d} \bs{v}, \quad
T(\bs{x}) = \frac{1}{d \rho(\bs{x})} \int_{\mathbb{R}^d} |\bs{v}-\bs{U}(\bs{x})|^2 f(\bs{x},\bs{v}) \,\mathrm{d} \bs{v}.
\end{equation}
Such a relaxation model eliminates the high-dimensional integral in $\mathcal{Q}[f,f]$, leading to a significantly lower computational cost. More accurate relaxation models include the ES-BGK model \cite{Holway1966} and the Shakhov model \cite{Shakhov1968}, which can better match the Prandtl number of the gas that is predicted wrongly in the BGK model. To solve the steady-state Boltzmann equation, one possibility is to consider the time-dependent Boltzmann equation and evolve the solution with time until the solution hardly changes. The classical DSMC (direct simulation of Monte Carlo) method \cite{Bird1994} can be considered this type of approach. Other numerical schemes solving the time-dependent Boltzmann equation include the work in \cite{Mieussens2000}, which considers the BGK collision operator and allows taking large time steps and applying implicit schemes to achieve fast convergence towards the steady-state solution. The UGKS (Unified gas-kinetic scheme) and UGKWP (Unified gas-kinetic wave-particle) method are also efficient schemes providing accurate approximations of numerical fluxes \cite{Xu2010, Liu2020}. More recently, some low-rank methods have been applied to tackle the high dimensionality of the Boltzmann equation \cite{Einkemmer2021, Hu2022}.

Some iterative methods solving the Boltzmann equation directly have also been developed. For instance, in \cite{Hu2019}, the Boltzmann equation is discretized using the moment method and solved using the nonlinear multi-level method. In \cite{Su2020Can}, an efficient iterative solver is proposed for the linear Boltzmann equation to achieve uniform convergence rate in all regimes, and the method has been generalized to the nonlinear Boltzmann equation \cite{Zhu2021} and polyatomic cases \cite{Zeng2023}. Such an approach is known as the general synthetic iterative scheme (GSIS), which takes the idea of the synthetic method for the linear radiative transfer equation \cite{Bond1983} and uses the diffusion limit of the Boltzmann equation as a preconditioner of the iterative method.

When using GSIS to solve the nonlinear Boltzmann equation, we are required to solve the nonlinear Navier-Stokes-Fourier (NSF) equation to obtain accurate approximations of conservative variables. The NSF equations involves the spatial variable $\bs{x}$ only, and are therefore cheaper to solve. However, in practice, solving the steady-state NSF equations is not an easy task at all and intricate codes may need to be written to achieve good convergence \cite{Rana2013}. The aim of our work is to develop a numerical solver that does not require solving Euler or NSF equations but can still achieve a convergence rate that is independent of $\epsilon$. Our numerical method is based on the symmetric Gauss-Seidel method applied to the spatial grid, and when updating the distribution function on each grid cell, a preconditioner that is based on the asymptotic limit and only requires solving a nonlinear algebraic system is proposed to accelerate the convergence of the inner iteration. The method is applied to both the BGK model and the binary collision to show its effectiveness. Furthermore, we couple the symmetric Gauss-Seidel method with the nonlinear multigrid technique to accelerate the convergence.

Below we will introduce our discretization of the Boltzmann equation in \Cref{sec:discretization}, and our numerical method combining the symmetric Gauss-Seidel method and the preconditioned fixed-point iteration will be introduced in \Cref{sec:numerical}. In \Cref{sec:experiments}, we will present our numerical experiments showing the effectiveness of our method. Some concluding remarks will be given in \Cref{sec:conclusion}.

\section{Discretization of the steady-state Boltzmann equation}
\label{sec:discretization}
In this section, we will demonstrate our discretization of the Boltzmann equation. In general, we will apply the discrete velocity model in the discretization of $\bs{v}$, and the finite volume method will be used for spatial discretization. Details are given in the following subsections.

\subsection{Spatial discretization}
Due to the conservation properties of the Boltzmann equation, the finite volume method and the discontinuous-Galerkin method are commonly used in the spatial discretization \cite{Jaiswal2019, Cai2022}. Here we adopt the finite volume method with piecewise linear construction, which can achieve second order of accuracy. For simplicity, below we introduce our discretization based on the spatially one-dimensional case, which can be extended to the two-dimensional case with rectangular grids without much difficulty.

Assume that $\bs{x} = (x_1, \cdots, x_d)^T$ and $\bs{v} = (v_1, \cdots, v_d)^T$, and the distribution function $f$ is homogeneous in $x_2, \ldots, x_d$. The Boltzmann equation \eqref{eq:Boltz} can be simplified to
\begin{equation} \label{eq:Boltz1D}
v_1 \frac{\partial}{\partial x} f(x,\bs{v}) = \frac{1}{\epsilon} \mathcal{Q}[f,f](x,\bs{v}), \qquad x \in \Omega \subset \mathbb{R}, \quad \bs{v} \in \mathbb{R}^d,
\end{equation}
where $x$ refers to $x_1$. Suppose the spatial domain $\Omega = (x_L, x_R)$ is divided into $N$ uniform grid cells of size $\Delta x$, where the $j$th cell is denoted by $[x_{j-1/2}, x_{j+1/2}]$. By the finite volume method, we use $\bar{f}(\bs{v})$ to denote the cell average of $f(x,\bs{v})$. The upwind method can then be applied to discretize the Boltzmann equation:
\begin{equation} \label{eq:upwind}
v_1^+ \frac{f^-_{j+1/2}(\bs{v}) - f^-_{j-1/2}(\bs{v})}{\Delta x} + v_1^- \frac{f^+_{j+1/2}(\bs{v}) - f^+_{j-1/2}(\bs{v})}{\Delta x} 
= \frac{1}{\epsilon} \mathcal{Q}[\bar{f}_j, \bar{f}_j].
\end{equation}
Here $v_1^+ = \max \{v_1,0\}$, $v_1^- = \min \{v_1,0 \}$, and the right-hand side should be replaced by $\epsilon^{-1} \mathcal{Q}^{\mathrm{BGK}}[\bar{f}_j]$ if the BGK collision term is applied. In the semi-discrete equation \eqref{eq:upwind}, the functions $f_{j+1/2}^{\pm}$ denotes the left and right limits of the distribution function at $x_{j+1/2}$ after linear reconstruction:
\begin{equation} \label{eq:reconstruction}
f^-_{j+1/2}(\bs{v}) = \bar{f}_j(\bs{v}) + \frac{1}{2} s_j(\bs{v}), \qquad
f^+_{j+1/2} = \bar{f}_{j+1}(\bs{v}) - \frac{1}{2} s_{j+1}(\bs{v}).
\end{equation}
For first-order schemes, we can simply choose $s_j(\bs{v}) \equiv 0$ for all $j = 1,\ldots,N$. To achieve second order of accuracy, the following choice of $s_j$ will be applied: 
\begin{equation} \label{eq:slope}
s_{j} = \left \{
\begin{aligned}
& \bar{f}_{j+1} - \bar{f}_j, \quad j=1; \\
& \frac{\bar{f}_{j+1}-\bar{f}_{j-1}}{2}, \quad j=2,\ldots,N-1; \\
& \bar{f}_j - \bar{f}_{j-1}, \quad j=N.
\end{aligned}
\right.
\end{equation}
Note that no limiters are applied in our current implementation. It will be shown that our method also works when limiters are applied.

\subsection{Velocity discretization} \label{sec:velocity_discretization}
To discretize the velocity variable, we pick a collection of points in $\mathbb{R}^d$, denoted by $\bs{v}_1, \ldots, \bs{v}_K$, and we assume that these points can be used as quadrature points, and the corresponding weights are $w_1, \ldots, w_K$ such that
\begin{equation} \label{eq:velo_integral}
\int_{\mathbb{R}^d} \varphi(\bs{v}) \,\mathrm{d}\bs{v} \approx \sum_{k=1}^K w_k \varphi(\bs{v}_k), \qquad \forall \varphi \in L^1(\mathbb{R}^d).
\end{equation}
Thus, the unknowns in the fully discrete equations are $\bar{f}_{j,k}$, which approximates $\bar{f}_j(\bs{v}_k)$, and \eqref{eq:upwind} is approximated by
\begin{equation} \label{eq:fully_discrete}
v_{k,1}^+ \frac{f_{j+1/2,k}^- - f_{j-1/2,k}^-}{\Delta x} + v_{k,1}^- \frac{f_{j+1/2,k}^+ - f_{j-1/2,k}^+}{\Delta x} = \frac{1}{\epsilon} \mathcal{Q}_k[\bar{\bs{f}}_j, \bar{\bs{f}}_j],
\end{equation}
where $\bar{\bs{f}}_j$ denotes the vector $(f_{j,1}, \ldots, f_{j,K})^T$, and the definitions of the reconstructed boundary values $f_{j+1/2,k}^{\pm}$ are similar to those of $f_{j+1/2}^{\pm}(\bs{v})$ as in \eqref{eq:reconstruction} and \eqref{eq:slope}.

For the BGK collision term, the right-hand side of \eqref{eq:fully_discrete} should be changed to
\begin{equation} \label{eq:BGK_discrete}
\frac{1}{\epsilon} \mathcal{Q}^{\mathrm{BGK}}_k[\bar{\bs{f}}_j] := \frac{\nu_j}{\epsilon} (\mathcal{M}_k[\bar{\bs{f}}_j] - f_{j,k}).
\end{equation}
In the BGK collision term \eqref{eq:BGK_collis}, the Maxwellian $\mathcal{M}[f]$ satisfies
\begin{displaymath}
\int_{\mathbb{R}^d} \bs{\phi}(\bs{v}) \mathcal{M}[f](\bs{v}) \,\mathrm{d}\bs{v} =
\int_{\mathbb{R}^d} \bs{\phi}(\bs{v}) f(\bs{v}) \,\mathrm{d}\bs{v},
\end{displaymath}
where $\bs{\phi}(\bs{v}) = (1,v_1,\ldots,v_d,|\bs{v}|^2)^T$, so that the collision does not change the local density, momentum and energy. Similarly, in the discrete form \eqref{eq:BGK_discrete}, we require $\mathcal{M}_k[\bar{\bs{f}}_j]$ to be a discrete Gaussian satisfying
\begin{equation} \label{eq:discrete_conservation}
\sum_{k=1}^K w_k \bs{\phi}(\bs{v}_k) \mathcal{M}_k[\bar{\bs{f}}_j] =
\sum_{k=1}^K w_k \bs{\phi}(\bs{v}_k) f_{j,k},
\end{equation}
where $\mathcal{M}_k[\bar{\bs{f}}_j]$ has the general form
\begin{equation} \label{eq:discrete_Gaussian}
\mathcal{M}_k[\bar{\bs{f}}_j] = \exp\left( \alpha_j + \bs{\beta}_j \cdot \bs{v}_k - \gamma_j |\bs{v}_k|^2 \right),
\end{equation}
in which the coefficients $\alpha_j$, $\bs{\beta}_j$ and $\gamma_j$ can be determined by plugging \eqref{eq:discrete_Gaussian} into \eqref{eq:discrete_conservation} and solving the resulting nonlinear system. Such a technique was introduced in \cite{Mieussens2000}. It is guaranteed that the discrete Maxwellian minimizes the discrete entropy of all distribution functions under the constraints \eqref{eq:discrete_conservation}.  In our implementation, the velocity domain $\mathbb{R}^d$ is truncated to $[-L,L]^d$, and a uniform grid applied to discretize the truncated domain, so that all $w_k$'s are equal to the volume of one grid cell. Newton's method is employed to solve the nonlinear system 
\eqref{eq:discrete_conservation}. The convergence can be achieved within five iterations on average.

For the quadratic collision term, the velocity domain is also truncated, and for the sake of efficiency, we also assume that the velocity domain is periodic and choose the set of discrete velocities to be
\begin{displaymath}
\Big\{ (k_1 \Delta v, \ldots, k_d \Delta v)^T \,\Big\vert\, k_1, \ldots, k_d \in \{-K/2+1, \cdots, K/2\} \Big\}.
\end{displaymath}
Thus, the Fourier spectral method can be applied to evaluate the collision term efficiently. We refer the readers to \cite{Mouhot2006, Gamba2017} for more details. In our work, the right-hand side of \eqref{eq:fully_discrete} is implemented as
\begin{equation} \label{eq:discrete_Q}
\mathcal{Q}_k[\bar{\bs{f}}_j, \bar{\bs{f}}_j] := \mathcal{Q}^{\mathrm{FSM}}_k[\bar{\bs{f}}_j, \bar{\bs{f}}_j] - \mathcal{Q}^{\mathrm{FSM}}_k\big[\bs{\mathcal{M}}[\bar{\bs{f}}_j], \bs{\mathcal{M}}[\bar{\bs{f}}_j]\big],
\end{equation}
where $\mathcal{Q}^{\mathrm{FSM}}_k$ denotes the collision term evaluated by the Fourier spectral method. Such a method ensures that the discrete Maxwellian $\bs{\mathcal{M}}[\bs{f}_j]$ is still the local equilibrium, and the same method has been adopted by \cite{Filbet2015, Cai2022}.

\section{Iterative method for the steady-state Boltzmann equation}
\label{sec:numerical}
In this section, we will discuss the solver of the steady-state Boltzmann equation. For simplicity, below we will present our method as an iterative method for the semidisctete equation \eqref{eq:upwind}, and it is straightforward to apply the velocity discretization in Section \ref{sec:velocity_discretization} and turn it into a solver of the fully-discrete equation \eqref{eq:fully_discrete}. This will also be commented at the end of this section.

\subsection{Symmetric Gauss-Seidel method for the BGK equation}
For the BGK equation, a classical iterative method to find the steady-state solution is the source iteration, which can be written as
\begin{equation} \label{eq:SI}
v_1^+ \frac{f^{-,(n+1)}_{j+1/2}(\bs{v}) - f^{-,(n+1)}_{j-1/2}(\bs{v})}{\Delta x} + v_1^- \frac{f^{+,(n+1)}_{j+1/2}(\bs{v}) - f^{+,(n+1)}_{j-1/2}(\bs{v})}{\Delta x} 
= \frac{\bar{\nu}_j^{(n)}}{\epsilon} \left( \mathcal{M} [\bar{f}_j^{(n)}](\bs{v}) - \bar{f}_j^{(n+1)}(\bs{v}) \right).
\end{equation}
Here the superscript ``$(n)$'' or ``$(n+1)$'' is the index for iteration steps. In each iteration, a linear system needs to be solved to determine $f_{j}^{(n+1)}$. This method is easy to implement but can converge slowly when the Knudsen number $\epsilon$ is small. By Fourier analysis, it can be found that the error amplification factor is $1 - O(\epsilon^2)$,  which rapidly approaches $1$ when $\epsilon$ gets smaller \cite{Su2020Fast}. The general synthetic iterative scheme \cite{Su2020Can} can achieve a uniform convergence behavior for a wide range of $\epsilon$, at the cost of solving the Navier-Stokes equations in each iteration. In our approach, we adopt the symmetric Gauss-Seidel (SGS) method as the skeleton of our method. This method can achieve an error amplification factor $1 - O(\epsilon)$, which performs better than source iteration. It is shown in \cite{Dong2024} that such an iterative method can also achieve an upper bound of the error amplification factor that is less than 1 for all $\epsilon$.

The SGS method for the semidiscrete BGK equation \eqref{eq:upwind} (with right-hand side replaced by $\epsilon^{-1} \mathcal{Q}^{\mathrm{BGK}}[\bar{f}_j]$ includes a forward scan and a backward scan:
\begin{equation} \label{eq:SGS_BGK}
\left \{
\begin{aligned}
& \frac{1}{\Delta x} \left( -v_1^+ \bar{f}^{(n+1/2)}_{j-1} + |v_1| \bar{f}^{(n+1/2)}_j + v_1^- \bar{f}^{(n)}_{j+1} \right)
- \frac{1}{2 \Delta x} \left( v_1^+ s^{(n)}_{j-1} - v_1 s^{(n)}_j + v_1^- s^{(n)}_{j+1} \right) \\
& \hspace{250pt} = \frac{\bar{\nu}_j^{(n)}}{\epsilon} \left( \mathcal{M} [\bar{f}^{(n+1/2)}_j] - \bar{f}^{(n+1/2)}_j \right), \\
& \frac{1}{\Delta x} \left( -v_1^+ \bar{f}^{(n+1/2)}_{j-1} + |v_1| \bar{f}^{(n+1)}_j + v_1^- \bar{f}^{(n+1)}_{j+1} \right)
- \frac{1}{2 \Delta x} \left( v_1^+ s^{(n+1/2)}_{j-1} - v_1 s^{(n+1/2)}_j + v_1^- s^{(n+1/2)}_{j+1} \right) \\
& \hspace{250pt} = \frac{\bar{\nu}_j^{(n+1/2)}}{\epsilon} \left( \mathcal{M} [\bar{f}^{(n+1)}_j] - \bar{f}^{(n+1)}_j \right).
\end{aligned}
\right.
\end{equation}
Note that the slope terms including $s_j(\bs{v})$ are treated explicitly by using only information from the previous scan, which allows the scheme to work even if nonlinear limiters are applied.

Such a method requires us to solve a nonlinear equation of $\bar{f}_j$ every time a grid cell is visited. In general, this nonlinear equation has the form
\begin{equation} \label{eq:inner_eqn}
\frac{|v_1|}{\Delta x} g(\bs{v}) + r(\bs{v}) = \frac{\nu}{\epsilon} (\mathcal{M} [g] (\bs{v}) - g(\bs{v})),
\end{equation}
where the function $r(\bs{v})$ and the constant $\nu$ are given, whereas $g(\bs{v})$ is the unknown function. For instance, during the forward scan, when the $j$th grid cell is visited,  we need to solve \eqref{eq:inner_eqn} with
\begin{equation} \label{eq:r}
\begin{gathered}
g(\bs{v}) = \bar{f}^{(n+1/2)}_j(\bs{v}), \qquad 
\nu = \bar{\nu}_j^{(n)}, \\
r(\bs{v})  = \frac{1}{\Delta x} \left( -v_1^+ \bar{f}^{(n+1/2)}_{j-1}(\bs{v}) + v_1^- \bar{f}^{(n)}_{j+1}(\bs{v}) - \frac{1}{2} \left( v_1^+ s^{(n)}_{j-1}(\bs{v}) - v_1 s^{(n)}_j(\bs{v}) + v_1^- s^{(n)}_{j+1}(\bs{v}) \right) \right).
\end{gathered}
\end{equation}
Here $\bar{f}^{(n+1/2)}_{j-1}$ is treated as a given function since the $(j-1)$th cell has been visited when visiting the $j$th cell. For backward scans, it can also be observed that the equation has the form \eqref{eq:inner_eqn}.

Solving \eqref{eq:inner_eqn} again requires iterations since $\mathcal{M}[g]$ is a nonlinear operator. One naive approach is to apply the following fixed-point iteration:
\begin{equation} \label{eq:fix_point}
\frac{|v_1|}{\Delta x} g^{(k+1)}(\bs{v}) + r(\bs{v}) = \frac{\nu}{\epsilon} (\mathcal{M} [g^{(k)}] (\bs{v}) - g^{(k+1)}(\bs{v})),
\end{equation}
so that
\begin{displaymath}
g^{(k+1)}(\bs{v}) = \left( \frac{|v_1|}{\Delta x} + \frac{\nu}{\epsilon} \right)^{-1} \left( \frac{\nu}{\epsilon} \mathcal{M}[g^{(k)}](\bs{v}) - r(\bs{v}) \right).
\end{displaymath}
In the following subsection, we will introduce solver of \eqref{eq:inner_eqn} which does not slow down when $\epsilon$ approaches zero.

\subsection{A preconditioned fixed-point iteration for the nonlinear equation}
\label{sec:synthetic}
The key idea to accelerate the fixed-point iteration is to precondition the solver using the asymptotic limit of the equation. Such an approach has been applied to kinetic equations in \cite{Bond1983, Su2020Fast} as synthetic methods. Here in \eqref{eq:inner_eqn}, when $\epsilon$ is small, we have $g = \mathcal{M}[g] + O(\epsilon)$. Our goal is to capture the equilibrium $\mathcal{M}[g]$ more accurately using this property.

The Maxwellian $\mathcal{M}[f]$ is the distribution function that has the same density, momentum and energy with $f$ and possesses the maximum entropy. Therefore, for $\bs{\varphi}(\bs{v}) = (1,v_1,\ldots,v_d,|\bs{v}|^2)^T$, it holds that
\begin{equation} \label{eq:maxw_moment}
\int_{\mathbb{R}^d} 
\bs{\varphi}(\bs{v})
\mathcal{M}[f] (\bs{x},\bs{v}) \,\mathrm{d} \bs{v}
= \int_{\mathbb{R}^d} 
\bs{\varphi}(\bs{v})
f(\bs{x},\bs{v}) \,\mathrm{d} \bs{v}.
\end{equation}
Using this property, we multiply both sides of \eqref{eq:inner_eqn} by $\bs{\varphi}(\bs{v})$ and integrate with respect to $\bs{v}$ to get
\begin{equation} \label{eq:macro}
\int_{\mathbb{R}^d} \bs{\varphi}(\bs{v}) \left[ \frac{|v_1|}{\Delta x} g(\bs{v}) + r(\bs{v}) \right] \mathrm{d}\bs{v} = 0.
\end{equation}
The system above contains $d+2$ equations, which can be used to solve the equilibrium part of $g(\bs{v})$. Following this idea, we write the equation above as
\begin{displaymath}
\int_{\mathbb{R}^d} \bs{\varphi}(\bs{v})  \frac{|v_1|}{\Delta x} \mathcal{M}[g](\bs{v}) \mathrm{d}\bs{v} +
\int_{\mathbb{R}^d} \bs{\varphi}(\bs{v}) \left( \frac{|v_1|}{\Delta x} \Big[g(\bs{v}) - \mathcal{M}[g](\bs{v})\Big] + r(\bs{v}) \right) \mathrm{d}\bs{v} = 0.
\end{displaymath}
Solving $\mathcal{M}[g]$ requires the knowledge of the non-equilibrium part $g - \mathcal{M}[g]$, which can be approximated using the result of the previous iteration in an iterative method. Then the solution of $\mathcal{M}[g]$ can be used to determine the updated $g$ according to the fixed-point iteration. In detail, the iteration goes as follows:
\begin{itemize}
\item Solve $\mathcal{M}^{(k+1)}(\bs{v}) = \exp\left(\alpha^{(k+1)} + \bs{\beta}^{(k+1)}  \bs{v} - \gamma^{(k+1)} |\bs{v}|^2\right)$ from
\begin{equation} \label{eq:step1}
\int_{\mathbb{R}^d} \bs{\varphi}(\bs{v}) \frac{|v_1|}{\Delta x} \mathcal{M}^{(k+1)}(\bs{v}) \,\mathrm{d}\bs{v} +
\int_{\mathbb{R}^d} \bs{\varphi}(\bs{v}) \left( \frac{|v_1|}{\Delta x} \Big[g^{(k)}(\bs{v}) - \mathcal{M}[g^{(k)}](\bs{v})\Big] + r(\bs{v}) \right) \mathrm{d}\bs{v} = 0.
\end{equation}
\item Compute $g^{(k+1)}$ by
\begin{equation} \label{eq:step2}
g^{(k+1)}(\bs{v}) = \left( \frac{|v_1|}{\Delta x} + \frac{\nu}{\epsilon} \right)^{-1} \left( \frac{\nu}{\epsilon} \mathcal{M}^{(k+1)}(\bs{v}) - r(\bs{v}) \right).
\end{equation}
\end{itemize}
When $\epsilon$ is close to zero, the solution $g(\bs{v})$ is dominated by its equilibrium part $\mathcal{M}[g](\bs{v})$. As a result, the first step of the iteration \eqref{eq:step1} can already provide a good estimation of $\mathcal{M}[g]$, leading to fast convergence of the method. Compared with the fixed-point iteration \eqref{eq:fix_point}, it does not suffer from the slowdown when $\epsilon$ gets smaller, so that a much better performance can be expected. For large $\epsilon$ ($\epsilon \geqslant 1$), our experiments also show that our new method requires slightly lower numbers of iterations. Due to the extra step \eqref{eq:step1} which needs several Newton iterations, the overall computational time of the two methods is generally comparable.

However, in practice, the solution to \eqref{eq:step1} may not exist. Let
\begin{displaymath}
(s_0, s_1, \ldots, s_d, s_{d+1})^T =
\int_{\mathbb{R}^d} \bs{\varphi}(\bs{v}) \left( \frac{|v_1|}{\Delta x} \Big[\mathcal{M}[g^{(k)}](\bs{v}) - g^{(k)}(\bs{v})\Big] - r(\bs{v}) \right) \mathrm{d}\bs{v}.
\end{displaymath}
Then the solution exists only if
\begin{equation} \label{eq:existence}
s_0 > 0 \quad \text{and} \quad s_0 s_{d+1} > \sum_{i=1}^d s_i^2,
\end{equation}
since $\mathcal{M}[g](\bs{v}) > 0$ and according to \eqref{eq:step1},
\begin{displaymath}
s_0 = \int_{\mathbb{R}^d} \frac{|v_1|}{\Delta x} \mathcal{M}[g](\bs{v}) \,\mathrm{d} \bs{v}, \qquad
s_0 s_{d+1} - \sum_{i=1}^d s_i^2 = \sum_{i=1}^d \int_{\mathbb{R}^d} s_0 \left(v_i - \frac{s_i}{s_0} \right)^2 \frac{|v_1|}{\Delta x} \mathcal{M}[g](\bs{v})\,\mathrm{d} \bs{v}.
\end{displaymath}
In our numerical experiments, the violation of \eqref{eq:existence} is rare but occurs occasionally. When it happens, we resort to the fixed-point iteration \eqref{eq:fix_point} to find the solution.

\begin{remark}
An alternative method to solve \eqref{eq:inner_eqn} when \eqref{eq:existence} is violated is to introduce a relaxation factor $\tau$, and the first step \eqref{eq:step1} is revised to
\begin{displaymath}
(1+\tau)\int_{\mathbb{R}^d} \bs{\varphi}(\bs{v})\frac{|v_1|}{\Delta x} \mathcal{M}^{(k+1)}(\bs{v}) \,\mathrm{d}\bs{v} +
\int_{\mathbb{R}^d} \bs{\varphi}(\bs{v}) \left( \frac{|v_1|}{\Delta x} \Big[g^{(k)}(\bs{v})-(1+\tau) \mathcal{M}[g^{(k)}](\bs{v}) \Big] + r(\bs{v}) \right) \mathrm{d}\bs{v} = 0.
\end{displaymath}
When $\tau = 0$, the equation above is identical to \eqref{eq:step1}. When $\tau \rightarrow +\infty$, the source iteration is recovered. In practice, one can choose $\tau \in [1,+\infty)$ to guarantee that the right-hand side satisfies \eqref{eq:existence}.
\end{remark}

\begin{remark}
Similar to the generalized synthetic iterative scheme (GSIS) proposed in \cite{Su2020Can}, here we also use the asymptotic limit to accelerate the convergence of the fixed-point iteration. In GSIS, such a technique is applied to the Boltzmann equation whose asymptotic limit is Navier-Stokes-Fourier equations, while we first discretize the Boltzmann equation, and the acceleration technique is applied only to the local problem inside the SGS iterations. Two key differences can be noticed here: first, the equation \eqref{eq:step1} is a system of nonlinear algebraic equations, which is easier to solve compared with the second-order Navier-Stokes-Fourier equations; second, GSIS needs to consider the consistency between the discretizations of the kinetic equations and the macroscopic equations, including the treatment of boundary conditions \cite{Liu2024}, whereas our method does not require such considerations. However, our approach suffers from the slower convergence for finer grids caused by the SGS iteration. We will use the multigrid technique to overcome this difficulty in our numerical experiments (see \Cref{sec:experiments} for more details).
\end{remark}

\subsection{Numerical method for the Boltzmann equation with quadratic collision operators}
We now generalize the method to the steady-state Boltzmmann equation with quadratic collision operators using the method of ``BGK penalty'' \cite{Filbet2010}. The idea is to split the collision term $\mathcal{Q}[f,f]$ into two parts:
\begin{displaymath}
\mathcal{Q}[f,f] = \mathcal{P}[f] + \mathcal{Q}^{\mathrm{BGK}}[f],
\end{displaymath}
where
\begin{equation} \label{eq:Pf}
\mathcal{P}[f] = \mathcal{Q}[f,f] - \nu(\mathcal{M}[f] - f).
\end{equation}
The SGS method \eqref{eq:SGS_BGK} then turns out to be
\begin{equation} \label{eq:SGS}
\left \{
\begin{aligned}
& \frac{1}{\Delta x} \left( -v_1^+ \bar{f}^{(n+1/2)}_{j-1} + |v_1| \bar{f}^{(n+1/2)}_j + v_1^- \bar{f}^{(n)}_{j+1} \right)
- \frac{1}{2 \Delta x} \left( v_1^+ s^{(n)}_{j-1} - v_1 s^{(n)}_j + v_1^- s^{(n)}_{j+1} \right) \\
& \hspace{210pt} = \frac{1}{\epsilon} \mathcal{P}[\bar{f}_j^{(n)}] + \frac{\bar{\nu}_j^{(n)}}{\epsilon} \left( \mathcal{M} [\bar{f}^{(n+1/2)}_j] - \bar{f}^{(n+1/2)}_j \right), \\
& \frac{1}{\Delta x} \left( -v_1^+ \bar{f}^{(n+1/2)}_{j-1} + |v_1| \bar{f}^{(n+1)}_j + v_1^- \bar{f}^{(n+1)}_{j+1} \right)
- \frac{1}{2 \Delta x} \left( v_1^+ s^{(n+1/2)}_{j-1} - v_1 s^{(n+1/2)}_j + v_1^- s^{(n+1/2)}_{j+1} \right) \\
& \hspace{190pt} = \frac{1}{\epsilon} \mathcal{P}[\bar{f}_j^{(n+1/2)}] + \frac{\bar{\nu}_j^{(n+1/2)}}{\epsilon} \left( \mathcal{M} [\bar{f}^{(n+1)}_j] - \bar{f}^{(n+1)}_j \right).
\end{aligned}
\right.
\end{equation}
Note that we have treated the term $\mathcal{P}[\cdot]$ ``explicitly'' to avoid solving nonlinear equations involving the binary collision operator. When $\epsilon$ is small, the solution $f$ is close to the local Maxwellian $\mathcal{M}[f]$, so that $\mathcal{P}[f] / \epsilon$ has magnitude $O(1)$. In this case, the explicit treatment of $\mathcal{P}[f]$ can be understood as an analog of the ``relaxed SGS method'' which solves the linear system $(D-L-U)x = b$ using the iteration $x^{(n+1)} = (D+D'-U)^{-1} [(L+D') (D+D'-L)^{-1} (U x^{(n)} + D' x^{(n)} + b) + b]$ instead of the classical SGS method $x^{(n+1)} = (D-U)^{-1} [L (D-L)^{-1} (U x^{(n)} + b) + b]$, where $D$ and $D'$ are diagonal matrices and $L$/$U$ are strictly lower-/upper-triangular.
One can still expect that the convergence does not slow down when $\epsilon$ approaches zero.

During the iteration, each equation in \eqref{eq:SGS} still has the form \eqref{eq:inner_eqn}. Compared with \eqref{eq:r}, the expression of $r(\bs{v})$ will include an extra term $-\epsilon^{-1} \mathcal{P}[\bar{f}_j^{(n)}]$, but the method introduced in \Cref{sec:synthetic} still applies.

\begin{remark}
The velocity discretization introduced in \Cref{sec:velocity_discretization} can be applied without difficulty. Note that the discrete Maxwellian is calculated using \eqref{eq:discrete_conservation} and \eqref{eq:discrete_Gaussian}, which conserves mass, momentum, energy exactly, so that \eqref{eq:macro} also holds exactly after discretization. It is also important to discretize the binary collision operator using \eqref{eq:discrete_Q}, which has a term $\mathcal{Q}_k^{\mathrm{FSM}}\big[\bs{\mathcal{M}}[\bar{\bs{f}}_j], \bs{\mathcal{M}}[\bar{\bs{f}}_j]\big]$ approximating zero. This guarantees that $\mathcal{P}[f]$ is zero when $f$ is the Maxwellian, so that the discrete $\mathcal{P}[f]/\epsilon$ indeed has magnitude $O(1)$ for small $\epsilon$.
\end{remark}

\begin{remark}
To further improve the convergence at small values of $\epsilon$, one can consider replacing the BGK collision term with the ES-BGK or Shakhov collision term. Using these models, one can choose $\nu$ appropriately so that $\mathcal{P}[f] \sim O(\epsilon^2)$, so that the extra term $\mathcal{P}[f]/\epsilon$ in \eqref{eq:SGS} approaches zero when $\epsilon$ tends to zero. This will be studied in our future work.
\end{remark}

\section{Numerical experiments}
\label{sec:experiments}
We are going to carry out some numerical experiments in this section to present the performances of the proposed methods for various Knudsen numbers $\epsilon$, especially when $\epsilon$ is small. Our experiments will involve the following methods for comparison:
\begin{itemize}
\item The source iteration \eqref{eq:SI}, abbreviated as SI below;
\item The SGS method with fixed-point iteration \eqref{eq:fix_point}, abbreviated as SGS-FP below;
\item The SGS method with preconditioned fixed-point iteration \eqref{eq:step1}\eqref{eq:step2}, abbreviated as SGS-PFP below;
\item The SGS-PFP method with multigrid acceleration, abbreviated as MG-SGS-PFP below.
\end{itemize}
Note that we have also implemented the multigrid method since the SGS iteration slows down significantly when the number of grid cells gets larger. Our implementation of the multigrid method follows the standard procedure of the nonlinear multigrid method \cite{Hackbusch1985}. Suppose $f_h$ is the numerical solution on the fine grid of cell size $h$ and the corresponding equation is
\begin{displaymath}
R_h(f_h) = r_h,
\end{displaymath}
where $R_h$ is a nonlinear operator coming from the discretization of the nonlinear equation, and $r_h$ includes contributions from boundary conditions.
To solve this problem, we first apply the SGS-PFP iteration a few times as the smoother to get an approximate solution $\bar{f}_h$. Afterwards, we solve the following coarse-grid equation of $f_H$:
\begin{equation} \label{eq:coarse}
R_H(f_H) = R_H(I_h^H \bar{f}_h) + I_h^H(r_h - R_h(\bar{f}_h)),
\end{equation}
where $R_H$ is the discrete nonlinear operator on the coarse grid, and $I_h^H$ is the restriction operator from the fine grid to the coarse grid. Once $f_H$ is solved, we correct $\bar{f}_h$ by
\begin{displaymath}
\hat{f}_h := \bar{f}_h + I_H^h(f_H - I_h^H \bar{f}_h),
\end{displaymath}
which completes one iteration. Note that $I_H^h$ stands for the prolongation operator from the coarse grid to the fine grid. The coarse-grid equation \eqref{eq:coarse} can also be solved by using an even coarser grid, resulting in a V-cycle multigrid method. In our implementation, we simply choose $H = 2h$. Some other parameters of the multigrid method are listed in \Cref{tab:mg}. In general, more smoothing steps are needed for second-order schemes according to our tests.
\begin{table}[!ht]
\centering
\caption{Parameters for the multigrid method.}
\label{tab:mg}
\begin{tabular}{ccccc}
\hline \hline
& Method & Coarsest grid & Pre-smoothing steps & Post-smoothing steps \\
\hline
\multirow{2}*{1D tests} & 1st-order scheme & $4$ cells & $1$ & $1$ \\
& 2nd-order scheme & $8$ cells & $5$ & $1$ \\
\hline
\multirow{2}*{2D tests} & 1st-order scheme & $5 \times 5$ cells & $1$ & $1$ \\
& 2nd-order scheme & $5 \times 5$ cells & $5$ & $1$ \\
\hline \hline
\end{tabular}
\end{table}

In all our experiments, we will consider problems with wall boundary conditions, and here we adopt the fully diffusive model defined by at a boundary point $\bs{x}^w$:
\begin{equation} \label{eq:boundary}
f(\bs{x}^w, \bs{v}) = \frac{\rho^w(\bs{x}^w)}{[2 \pi T^w(\bs{x}^w)]^{d/2}} \exp \left( - \frac{| \bs{v} - \bs{U}^w(\bs{x}^w) |^2}{2 T^w(\bs{x}^w)} \right), \qquad [\bs{v} - \bs{U}^w(\bs{x}^w)] \cdot \bs{n}(\bs{x}^w) < 0.
\end{equation}
In the equation above, $\bs{n}(\bs{x}^w)$ refers to the outer normal vector at $\bs{x}^w$, and $T^w(\bs{x}^w)$ and $\bs{U}^w(\bs{x}^w)$ denote the wall temperature and velocity, respectively. The density of the reflected particles $\rho^w(\bs{x}^w)$ is determined by the ``no mass flux'' condition on the wall:
\begin{displaymath}
\int_{\mathbb{R}^d} \Big([\bs{v} - \bs{U}^w(\bs{x}^w)] \cdot \bs{n}(\bs{x}^w) \Big) f(\bs{x}^w, \bs{v}) \,\mathrm{d}\bs{v} = 0,
\end{displaymath}
which yields
\begin{equation*}
\rho^w(\bs{x}^w) = \sqrt{\frac{2 \pi}{T^w(\bs{x}^w)}} \int_{[\bs{v} - \bs{U}^w(\bs{x}^w)] \cdot \bs{n}(\bs{x}^w) > 0} \Big([\bs{v} - \bs{U}^w(\bs{x}^w)] \cdot \bs{n}(\bs{x}^w) \Big) f(\bs{x}^w, \bs{v}) \,\mathrm{d} \bs{v}.
\end{equation*}
In our experiments, such boundary conditions are applied to all $\bs{x}^w \in \partial \Omega$. As a result, an additional condition specifying the total mass is needed to uniquely determine the solution, which reads
\begin{equation} \label{eq:mass_conserv}
\int_{\Omega} \int_{\mathbb{R}^d} f (\bs{x}, \bs{v}) \,\mathrm{d} \bs{v} \,\mathrm{d} \bs{x} = C,
\end{equation}
where $C$ is a given constant. Numerically, since the SGS method does not preserve the total mass, we apply a uniform scaling for all $\bar{f}_j$ after each iteration to maintain the property \eqref{eq:mass_conserv}. 
Unless otherwise specified, the SGS iteration or the source iteration to solve \eqref{eq:SGS_BGK} or \eqref{eq:SGS} will be terminated if the numerical solution satisfies
\begin{displaymath}
\sqrt{\int_{\Omega} \int_{\mathbb{R}^d} \left(\bs{v} \cdot \nabla_{\bs{x}} f - \frac{1}{\epsilon} \mathcal{Q}[f,f] \right)^2 \,\mathrm{d}\bs{v} \,\mathrm{d}\bs{x}} < 10^{-5},
\end{displaymath}
where the spatial and velocity variables are discretized according to the description in \Cref{sec:discretization}.
For inner iterations including SGS-FP \eqref{eq:fix_point} and SGS-PFP \eqref{eq:step1}\eqref{eq:step2} to solve \eqref{eq:inner_eqn}, we terminate the iteration when the residual
\begin{displaymath}
    \left\| \frac{|v_1|}{\Delta x} g(\bs{v}) + r(\bs{v}) - \frac{\nu}{\epsilon} (\mathcal{M}[g](\bs{v}) - g(\bs{v}) \right\|_{L^2(\mathbb{R}^d)}
\end{displaymath}
reaches $10^{-8}$. 

In the following two subsections, we will present our test results for the BGK collision model and the binary collision model, respectively.

\subsection{BGK collision model}

We consider the BGK collision model \eqref{eq:BGK_collis} where $\mathcal{M}[f]$ is the local Maxwellian defined in \eqref{eq:maxwellian}. We allow the spatial variable $\bs{x}$ and the velocity variable $\bs{v}$ to have different dimensions. For example, when $\bs{x} \in \mathbb{R}^2$ and $\bs{v} \in \mathbb{R}^3$, we call this the ``2D3V'' case. Below we will start with the 1D1V case, followed by higher dimensional simulations.

\subsubsection{One-dimensional heat transfer} \label{sec:heat_trans_1D1V}

A toy model in which both the spatial and velocity variables are one-dimensional is formulated by
\begin{equation*}
v \frac{\partial}{\partial x} f(x,v) = \frac{1}{\epsilon} (\mathcal{M}[f](x,v) - f(x,v)), \qquad x \in (-1/2, 1/2) \subset \mathbb{R}, \quad v \in \mathbb{R},
\end{equation*}
and the boundary conditions are given by $T^w(-1/2) = 1$, $T^w(1/2) = 2$ and $U^w(-1/2) = U^w(1/2) = 0$ in \eqref{eq:boundary}. For the numerical simulations, the spatial domain is discretized on a uniform grid with $\Delta x = 1/256$, and the velocity domain is truncated to $[-L, L]$ with $L = 6$ and discretized uniformly with $\Delta v$ = $(2L)/50$. 

\begin{figure}[!ht] \label{fig:rho_T_1D1V}
    \centering
    \subfigure[Density]{
    \includegraphics[width=0.4\textwidth, trim=2 2 30 15, clip]{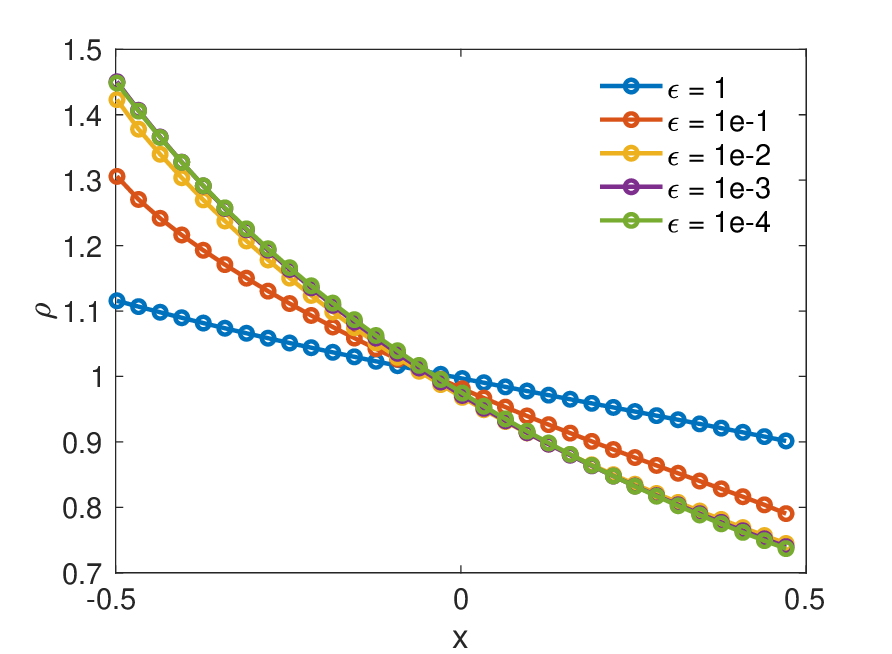}
    } \quad
    \subfigure[Temperature]{
    \includegraphics[width=0.4\textwidth, trim=2 2 30 15, clip]{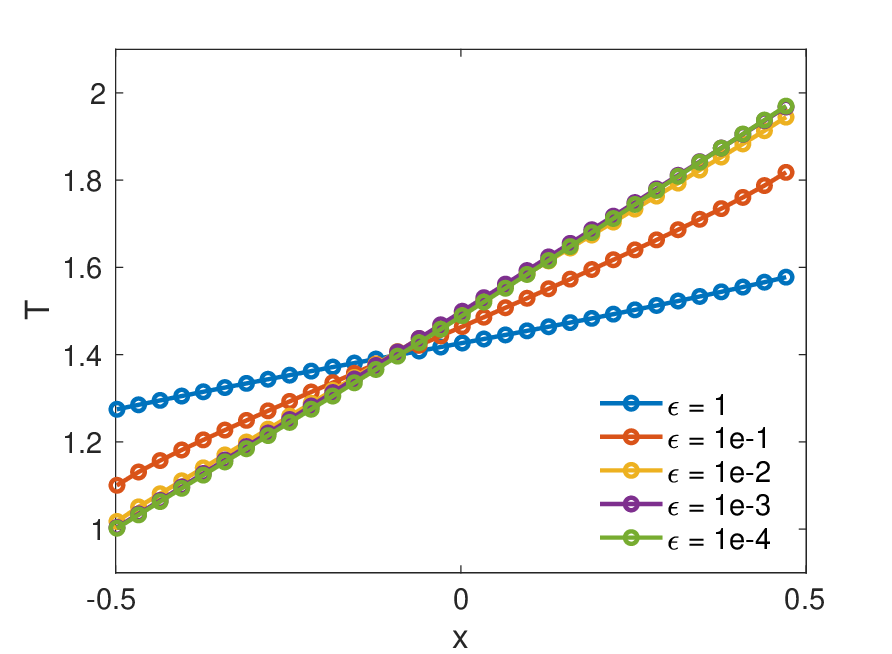}
    }
    \caption{The distributions of density and temperature for one-dimensional heat transfer.}
\end{figure}

The profiles of density $\rho(x)$ and temperature $T(x)$, which are defined in \eqref{eq:moment} and computed by the second-order scheme \eqref{eq:upwind}--\eqref{eq:slope}, are presented in \Cref{fig:rho_T_1D1V}. 
It is correctly predicted that when $\epsilon$ gets larger, the fluctuation of both density and temperature becomes weaker, and the temperature jump appears more significant. For $\epsilon = 10^{-3}$ and $10^{-4}$, the two curves in both figures almost coincide, as both results are already fairly close to the solutions of classical continuum fluid models. To test the order of convergence, the $L^2$ errors of density $\rho$ and temperature $T$ obtained from numerical solutions are evaluated:
\begin{equation*} 
\mathit{err} = \sqrt{\Delta x \sum_{j=1}^{N_x} \left( \bar{m}_j - \bar{m}^{\text{ref}}_j \right) } \end{equation*} 
where $m$ stands for $\rho$ or $T$, and $\bar{m}_j$ and $\bar{m}^{\text{ref}}_j$ are the numerical and reference cell averages of $m(x)$ in the $j$-th grid cell, respectively. Since the exact solution cannot be obtained analytically, we choose the numerical result with $\Delta x = 2^{-15}$ computed by the second-order scheme as the reference solution.
\begin{figure}[!ht] \label{fig:order}
    \centering
    \subfigure[1st-order scheme]{
    \includegraphics[width=0.405\textwidth, trim=2 2 30 15, clip]{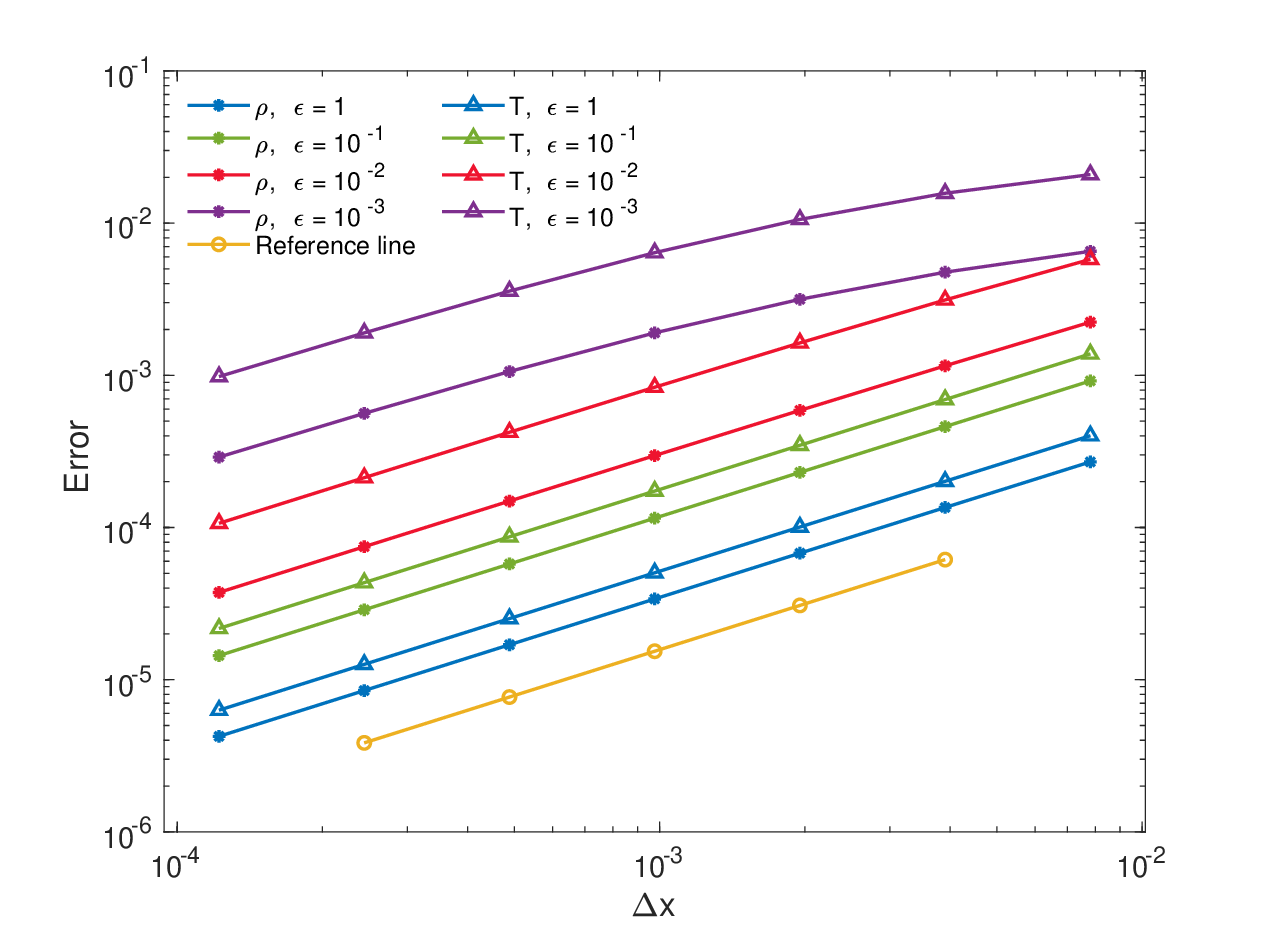}
    } \quad
    \subfigure[2nd-order scheme]{
    \includegraphics[width=0.4\textwidth, trim=2 2 30 15, clip]{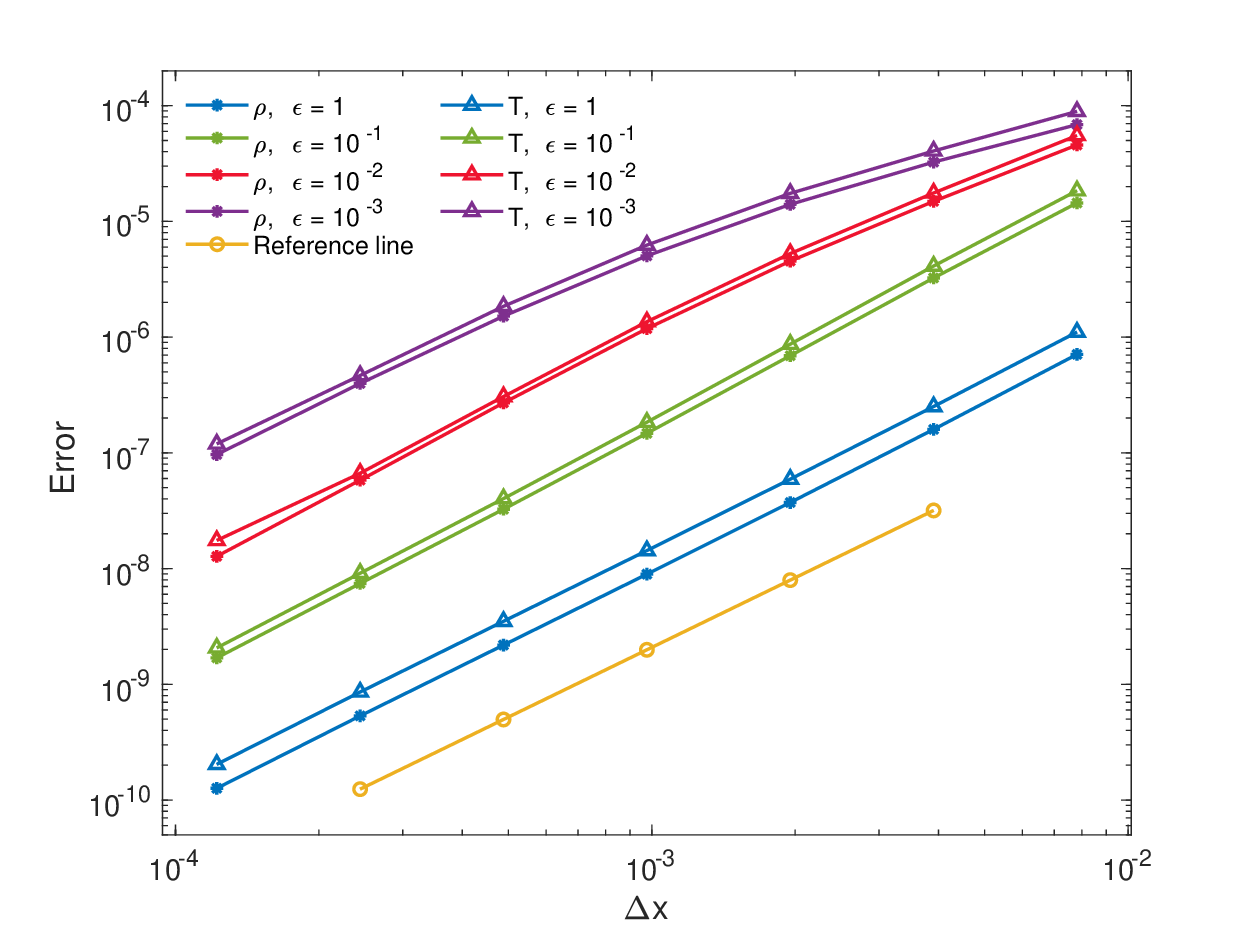}
    }
    \caption{The errors of upwind finite volume schemes with various $\epsilon$ for one-dimensional heat transfer.}
\end{figure}
The errors of first-order and second-order schemes are displayed in \Cref{fig:order} for $\Delta x = 1/(128 \cdot 2^j)$, $j = 0,1,\ldots,6$ and $\epsilon = 10^{-k}$, $k = 0,1,2,3$. In order to ensure sufficient accuracy, the iteration is terminated if the residual has a norm less than $10^{-9}$. We observe that the desired order of convergence can be achieved for various $\epsilon$.

\begin{figure}[!ht] \label{fig:rate_1D1V}
    \centering
    \subfigure[Source iteration (1st-order)]{\label{fig:rate_SI_1D1V}
    \includegraphics[width=0.31\textwidth, trim=2 2 35 20, clip]{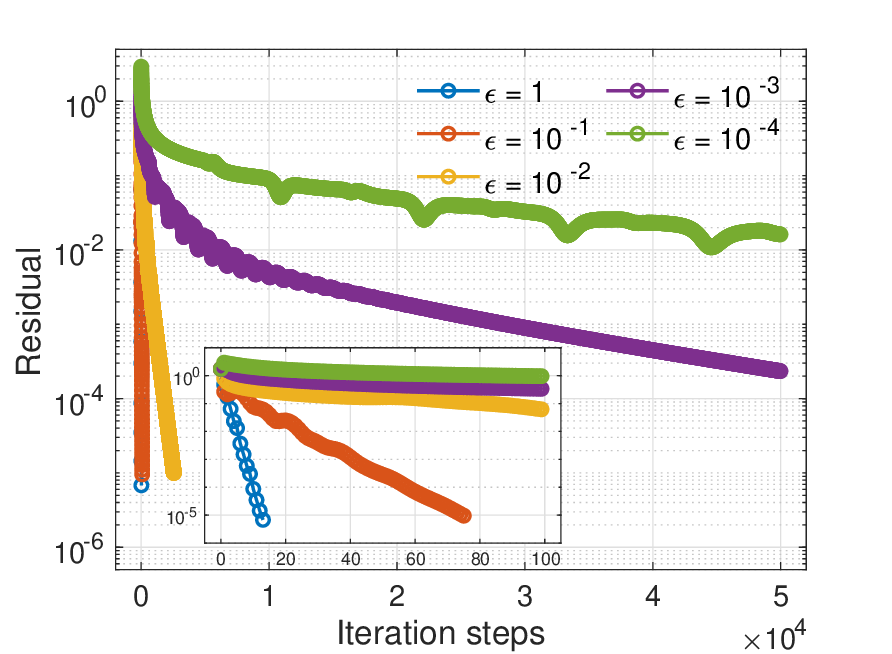}
    } 
    \subfigure[SGS method (1st-order)]{\label{fig:rate_SGS-1st_1D1V}
    \includegraphics[width=0.305\textwidth, trim=2 2 35 20, clip]{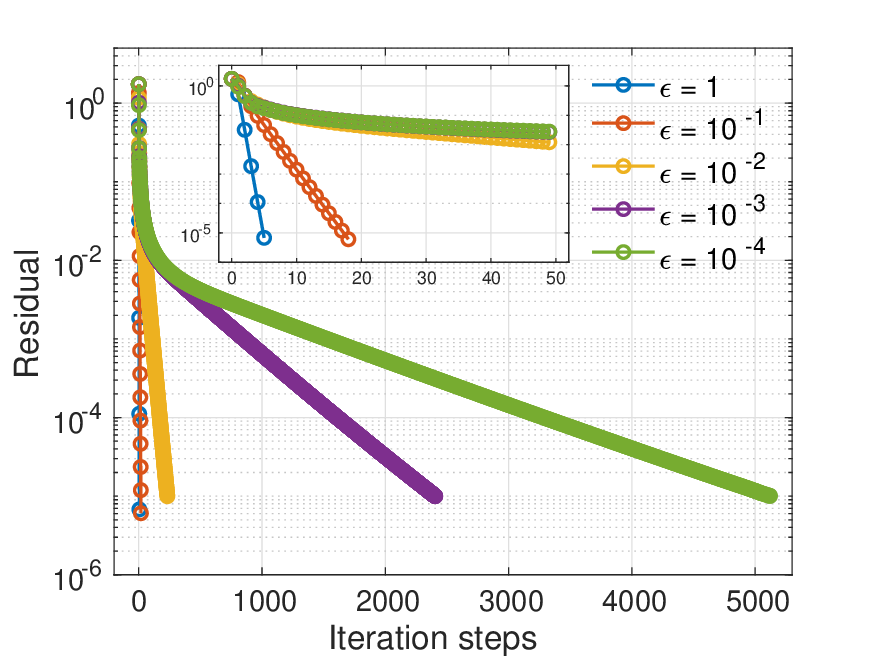}
    } 
    \subfigure[SGS method (2nd-order)]{\label{fig:rate_SGS-2nd_1D1V}
    \includegraphics[width=0.31\textwidth, trim=2 2 35 20, clip]{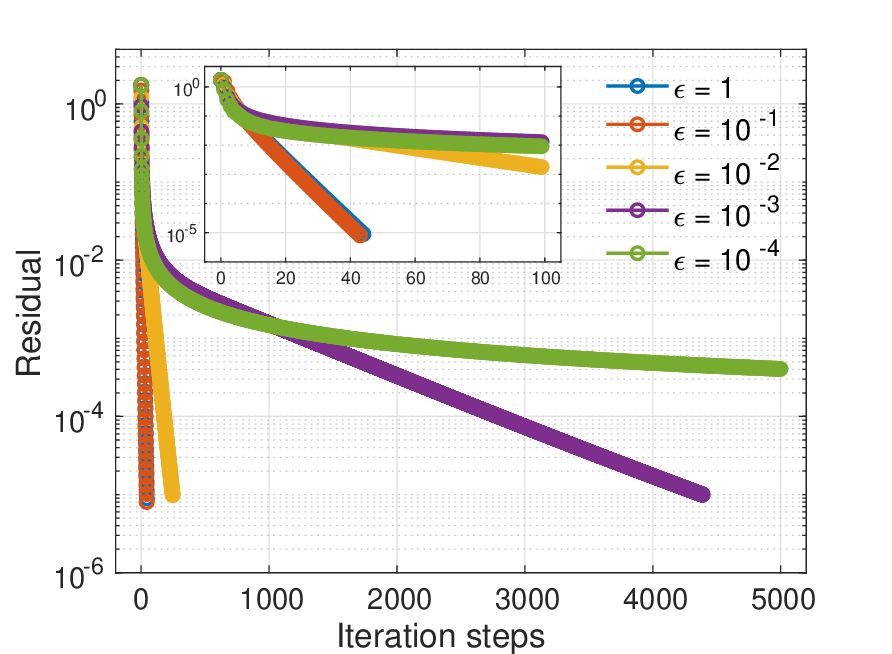}
    }
    \caption{The convergence of iterative methods for one-dimensional heat transfer.}
\end{figure}

In the following, we compare the efficiency of the source iteration and SGS method without using multigrid methods. We first compare the number of outer iterations, for which the decay of the residual norms is plotted in \Cref{fig:rate_1D1V}. Note that both SGS-FP and SGS-PFP lead to identical convergence plots since their difference is only the method for inner iterations.  For the first-order scheme, \Cref{fig:rate_SI_1D1V} and \ref{fig:rate_SGS-1st_1D1V} show that the number of iterations grows as $\epsilon$ decreases for both methods, but compared with the SGS method, the source iteration requires a significantly larger number of iterations, and the residual norm is not monotonically decreasing. Such a phenomenon is also observed in \cite{Dong2024}. The convergence of the second-order SGS method is slower due to the wider stencil in the numerical scheme (see \Cref{fig:rate_SGS-2nd_1D1V}), but the second-order method is also expected to achieve more accurate solutions. \Cref{fig:rate_SGS-1st_1D1V} and \Cref{fig:rate_SGS-2nd_1D1V} both show that the SGS method brings down the residual quickly in the first few steps, and then slows down gradually. This is the typical behavior of classical iterative methods, which eliminate high-frequency errors faster than low-frequency errors. We will later improve the performance by introducing multigrid methods.

\begin{figure}[!ht]  \label{fig:inner_1st_1D1V}
    \centering
    \subfigure[$\epsilon = 1$]{
    \label{fig:inner_1st_1D1V_0}
    \includegraphics[width=0.4\textwidth, trim=2 2 35 15, clip]{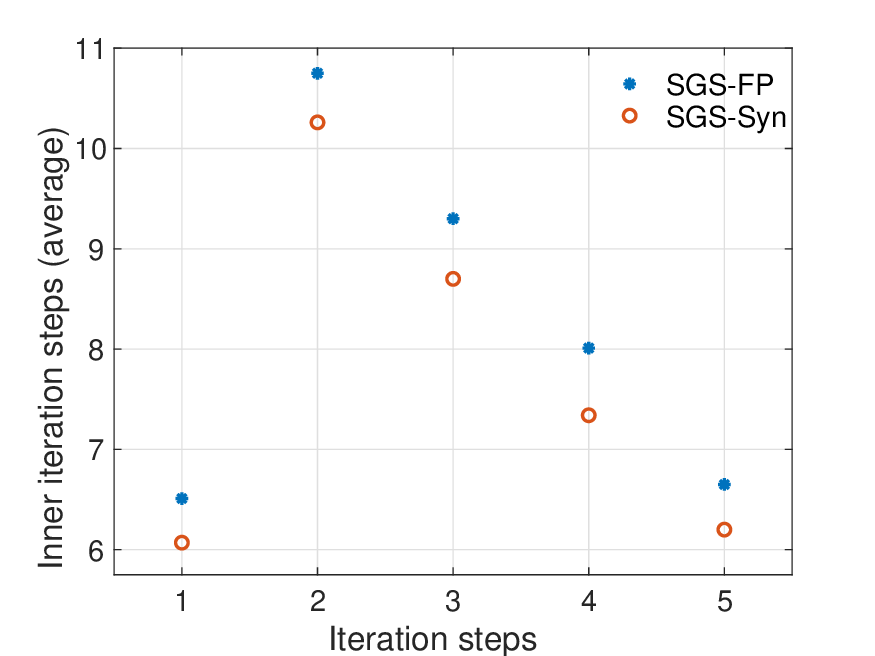}
    } \quad
    \subfigure[$\epsilon = 10^{-1}$]{
    \label{fig:inner_1st_1D1V_1}
    \includegraphics[width=0.4\textwidth, trim=2 2 35 15, clip]{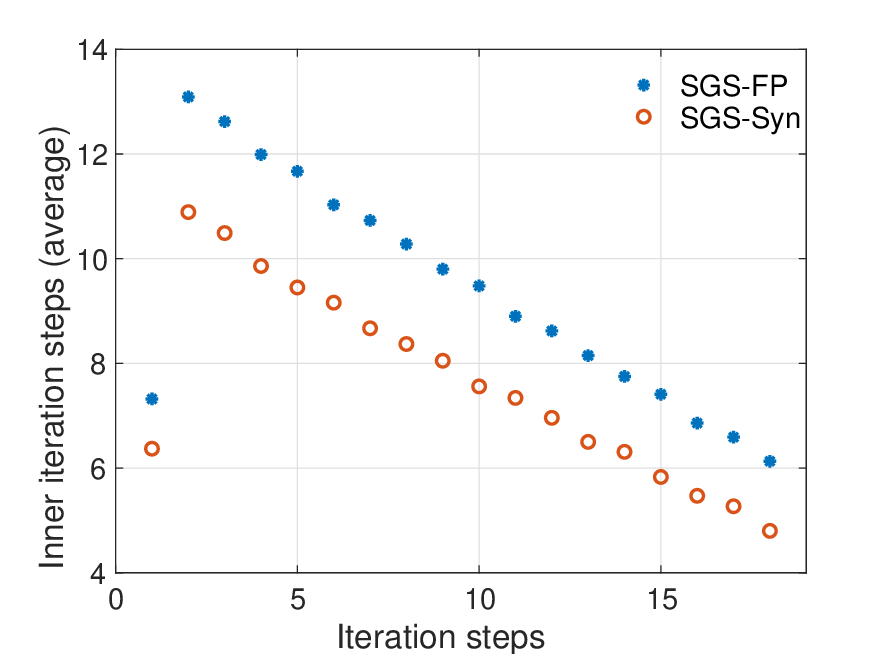}
    } \\
    \subfigure[$\epsilon = 10^{-2}$]{
    \includegraphics[width=0.4\textwidth, trim=2 2 35 15, clip]{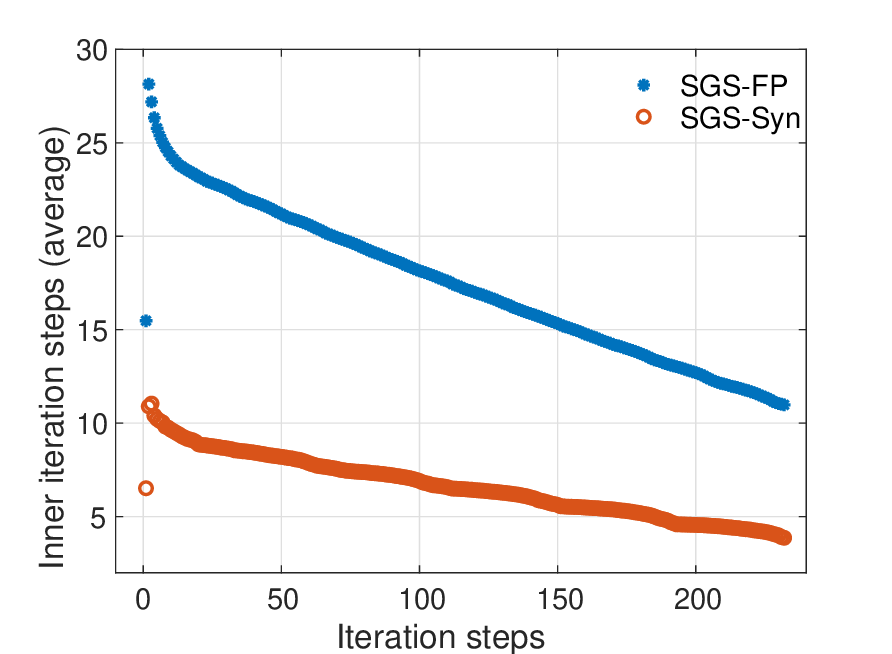}
    } \quad
    \subfigure[$\epsilon = 10^{-3}$]{
    \label{fig:inner_1st_1D1V_3}
    \includegraphics[width=0.4\textwidth, trim=2 2 35 15, clip]{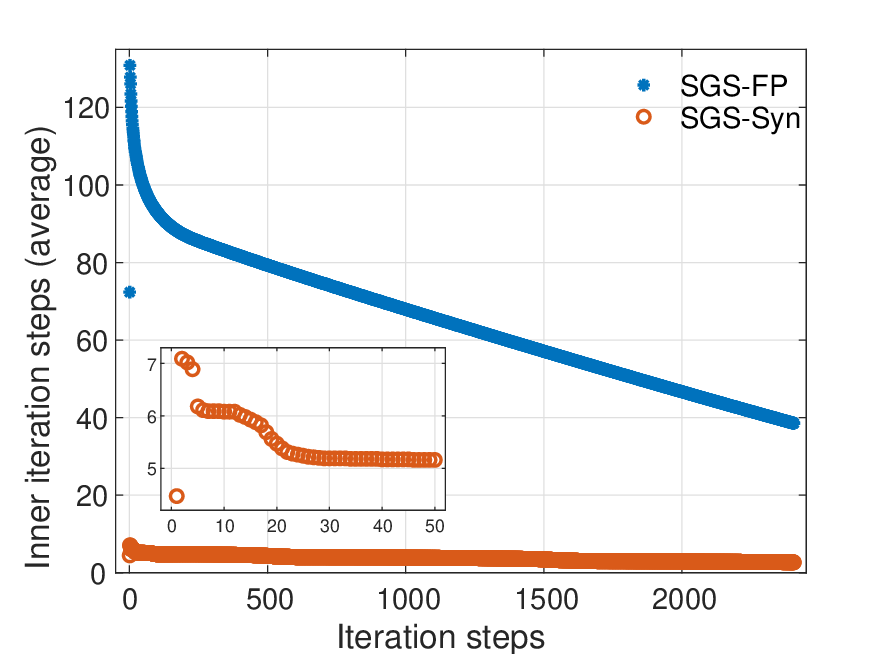}
    } \\
    \caption{The average inner iterations of iterative methods (1st-order) for one-dimensional heat transfer.}
\end{figure}

\begin{figure}[!ht] \label{fig:inner_2nd_1D1V}
    \centering
    \subfigure[$\epsilon = 1$]{
    \label{fig:inner_2nd_1D1V_0}
    \includegraphics[width=0.4\textwidth, trim=2 2 35 15, clip]{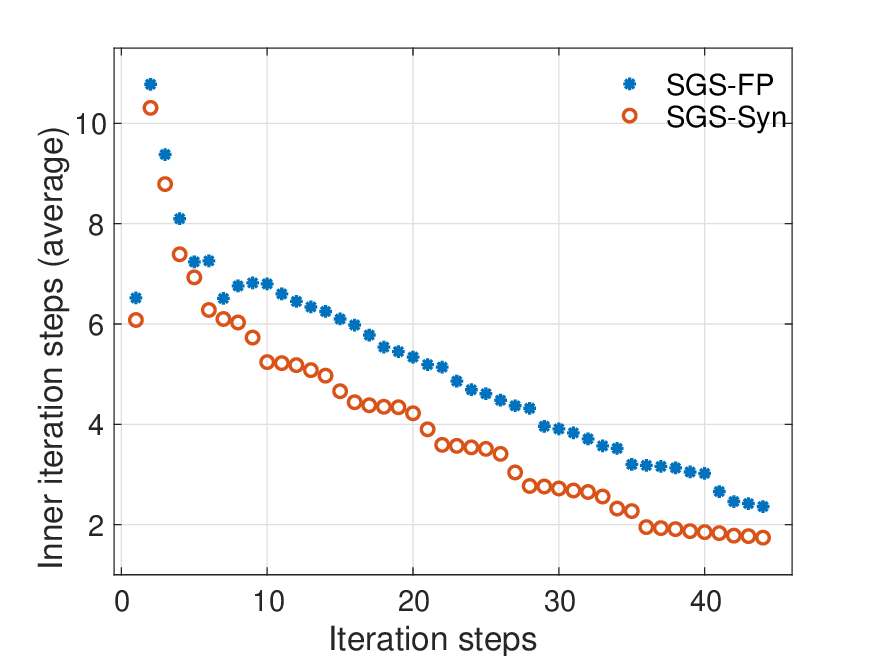}    
    } \quad
    \subfigure[$\epsilon = 10^{-1}$]{
    \includegraphics[width=0.4\textwidth, trim=2 2 35 15, clip]{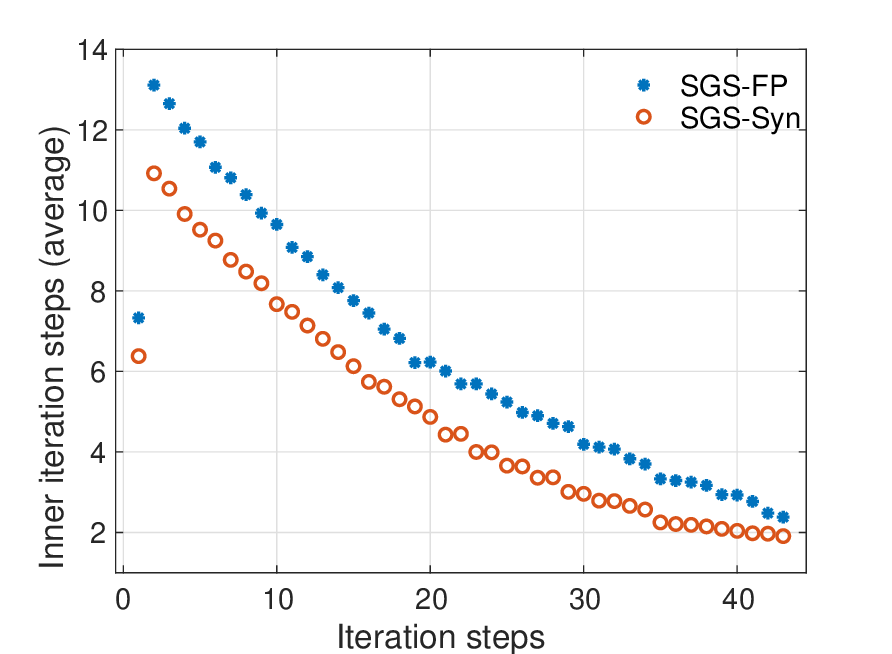}    
    } \\
    \subfigure[$\epsilon = 10^{-2}$]{
    \includegraphics[width=0.4\textwidth, trim=2 2 35 15, clip]{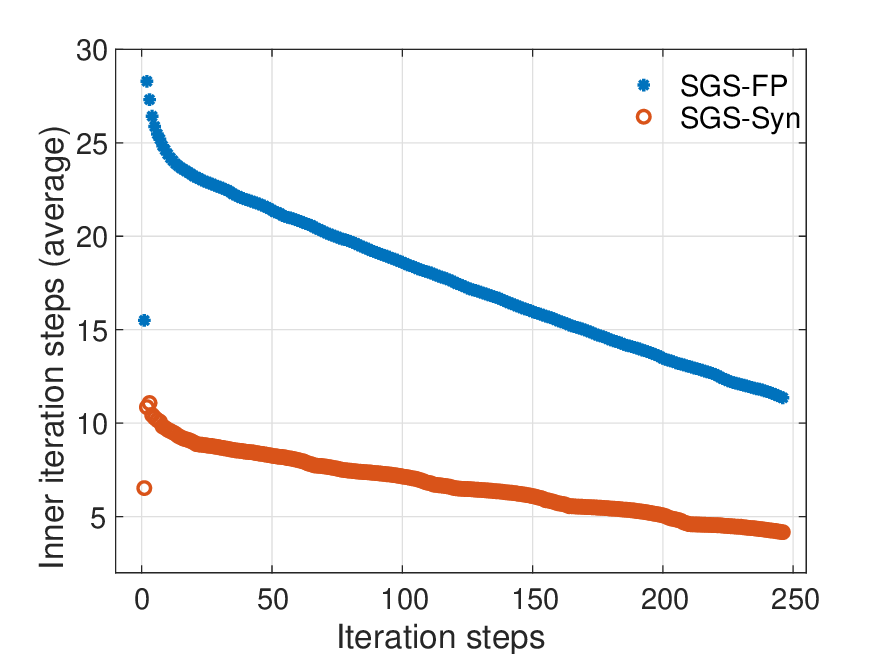}    
    } \quad
    \subfigure[$\epsilon = 10^{-3}$]{
    \label{fig:inner_2nd_1D1V_3}
    \includegraphics[width=0.4\textwidth, trim=2 2 35 15, clip]{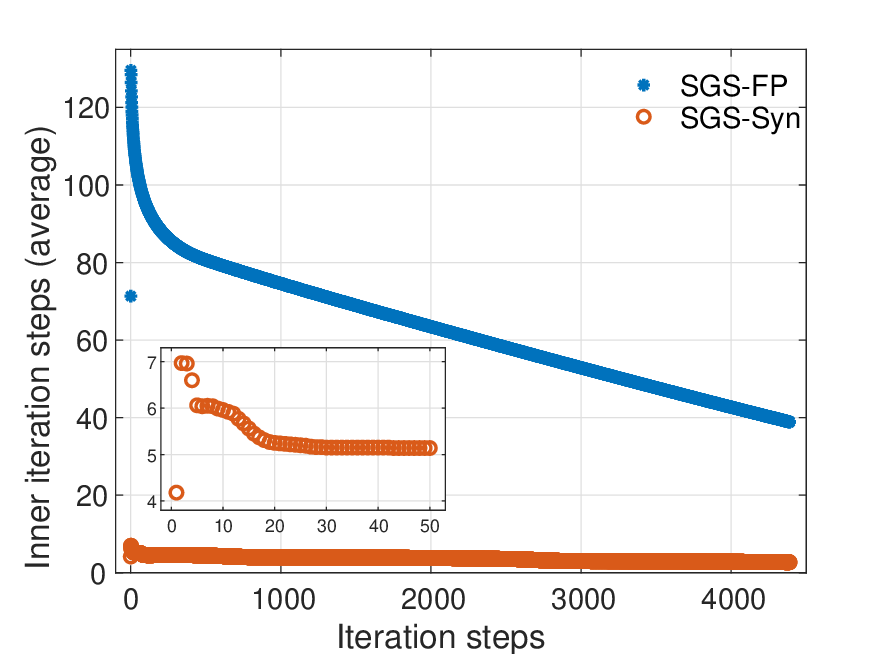}    
    } \\
    \caption{The average inner iterations of iterative methods (2nd-order) for one-dimensional heat transfer.}
\end{figure}

In order to illustrate how the preconiditioned method \eqref{eq:step1}\eqref{eq:step2} accelerates the simulation, we record the average number of inner iterations for both the SGS-FP method \eqref{eq:fix_point} and the SGS-PFP method in \Cref{fig:inner_1st_1D1V} and \ref{fig:inner_2nd_1D1V}, where the horizontal axis represents the index of the outer iteration (SGS iteration). Since the SGS method includes two sweeps over all grid cells, and during each sweep, every grid cell requires a few inner iterations to solve the nonlinear algebraic system, the vertical axis is actually the total number of inner iterations divided by twice the number of grid cells. It is shown that for large values of $\epsilon$ such as $1$ and $0.1$, the SGS-FP method and the SGS-PFP method have similar numbers of inner iterations for both first- and second-order schemes, and as the numerical solution gets closer to the exact solution, the number of inner iterations also reduces. When the value of $\epsilon$ reduces to $10^{-2}$, the SGS-FP method requires a lot more inner iterations, and it continues to grow as $\epsilon$ further decreases to $10^{-3}$. On the contrary, the SGS-PFP method, in both first- and second-order versions, shows a stable performance for all $\epsilon$, with the maximum number of inner iterations being less than $11$. In particular, when $\epsilon = 10^{-3}$, the maximum number of iterations is even less than $8$, since the local Maxwellian is dominant in the distribution function.

\begin{figure}[!ht] \label{fig:mg_1D1V}
    \centering
    \subfigure[1st-order scheme]{
    \includegraphics[width=0.4\textwidth, trim=2 2 40 20, clip]{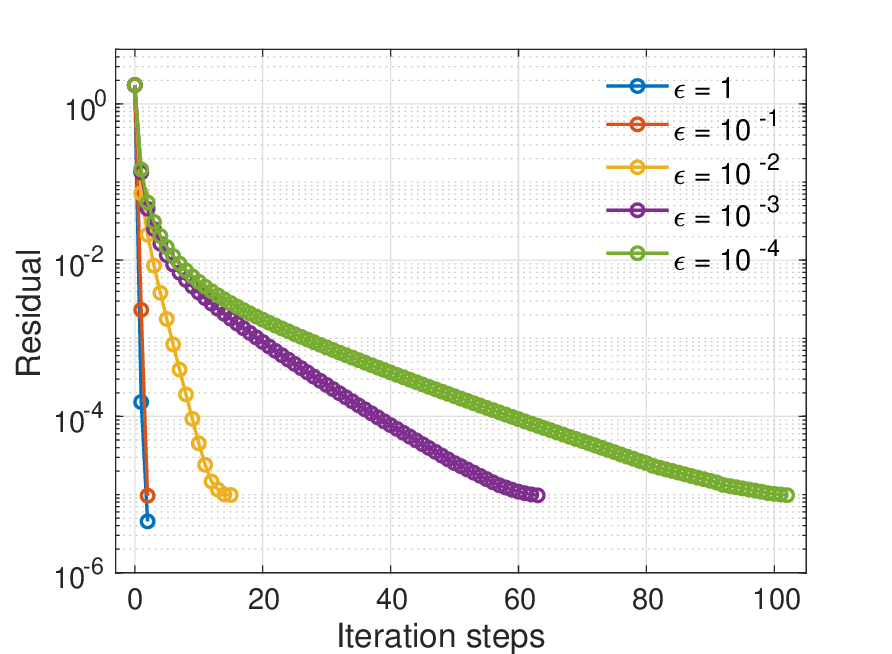}
    } \quad
    \subfigure[2nd-order scheme]{
    \includegraphics[width=0.4\textwidth, trim=2 2 40 20, clip]{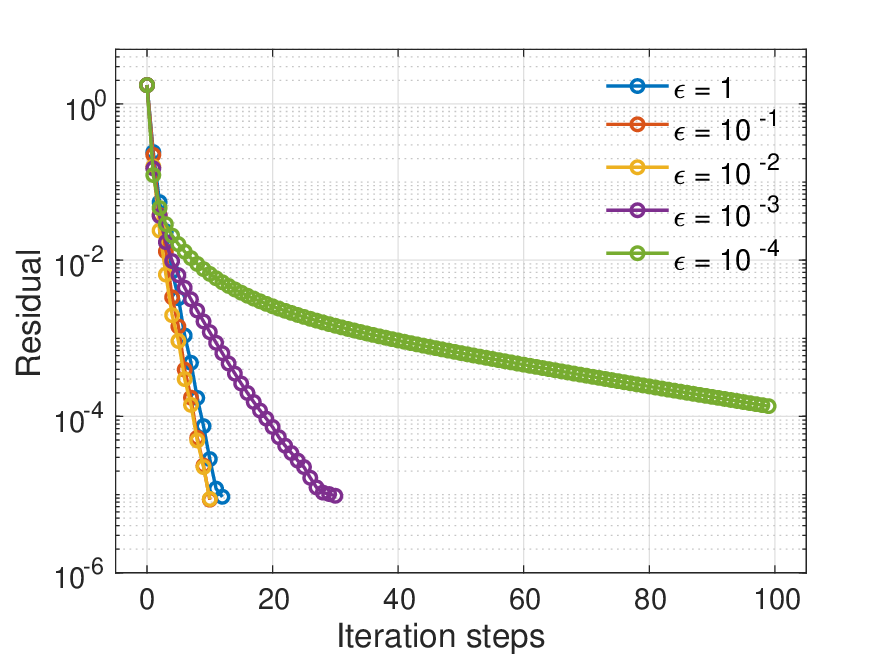}
    }
    \caption{The convergence of multigrid method for one-dimensional heat transfer.}
\end{figure}

In general, for large $\epsilon$, the transport term in the Boltzmann-BGK equation is dominant, and by sweeping, the steady-state can be quickly achieved. But when $\epsilon$ is small, by Chapman-Enskog expansion, it can be derived that the moments of the distribution function show some elliptic behavior, and thus the convergence of the SGS method may slow down significantly as the number of grid points increases. This explains why smaller $\epsilon$ requires more outer iterations. To alleviate this problem, we also test the multigrid technique introduced in the beginning of \Cref{sec:experiments}, and the results of the MG-SGS-PFP method are presented in \Cref{fig:mg_1D1V}. Although the number of iterations still increases as $\epsilon$ approaches zero, the convergence rate is remarkably improved.

Finally, to demonstrate the efficiency more intuitively for all the above methods, in \Cref{tab:time_1D1V}, we compare the actual computational time under the same computing environment.
\begin{table}[!ht] \label{tab:time_1D1V}
    \caption{Computing time for one-dimensional heat transfer (measure time by the second).}
    \footnotesize
    \centering
    \begin{tabular}{cccccccc}
        \hline 
        \hline 
        Order & & & $\epsilon = 1$ & $\epsilon = 10^{-1}$ & $\epsilon = 10^{-2}$ & $\epsilon = 10^{-3}$ & $\epsilon = 10^{-4}$ \\
        \hline
        \multirow{4}{*}{First} & SI & & $0.17$ & $0.96$ & $31.02$ & $1236.77$ & $18872.53$  \\
& SGS-FP & & $0.64$ & $2.73$ & $62.67$ & $2159.76$ & $27578.61$ \\
& SGS-PFP & & $0.95$ & $3.24$ & $35.79$ & $216.01$ & $276.23$ \\
& MG-SGS-PFP & & $0.84$ & $1.47$ & $8.72$ & $24.34$ & $26.35$ \\
        \hline
        \multirow{3}{*}{Second} & SGS-FP & & $3.91$ & $4.57$ & $66.73$ & $3974.22$ & $>36000$ \\
& SGS-PFP & & $4.54$ & $5.21$ & $39.80$ & $391.63$ & $1855.68$ \\
& MG-SGS-PFP & & $8.04$ & $8.15$ & $10.36$ & $27.91$ & $99.58$ \\
        \hline
        \hline 
    \end{tabular}
\end{table}
For first-order schemes with $\epsilon = 1$ and $0.1$, the classical SI has the best performance due to its simple implementation and least work on each grid cell. The computational time of SGS-FP is less than SGS-PFP since each inner iteration of SGS-FP (see \eqref{eq:fix_point}) is cheaper than that of SGS-PFP (see \eqref{eq:step1}\eqref{eq:step2}), despite a slightly larger average number of iterations in SGS-FP according to \Cref{fig:inner_1st_1D1V_0} and \ref{fig:inner_1st_1D1V_1}, and the multigrid technique does not provide much help. But when $\epsilon$ gets smaller, the SGS-PFP gets more advantageous. For example, when $\epsilon = 10^{-3}$, the computational time of SGS-PFP is only about $17\%$ of the time cost by SI, and this ratio further reduces to $1.5\%$ when $\epsilon = 10^{-4}$. The SGS-FP method has a worse performance than SI and should not be considered. Meanwhile, the multigrid method starts to take effect, which cuts down the computational time by around $90\%$. For second-order schemes, the SI was not tested, but the computational time of the other three methods follows the same trend. Note that for $\epsilon = 10^{-3}$ and $10^{-4}$, the multigrid technique is even more effective than in the first-order methods.

\subsubsection{Heat transfer in a cavity}  \label{sec:BGK_2D3V}

\begin{figure}[!ht] \label{fig:rho_2D3D_BGK}
    \centering
    \subfigure[$\epsilon = 1$]{
    \includegraphics[width=0.31\textwidth, trim=20 2 50 20, clip]{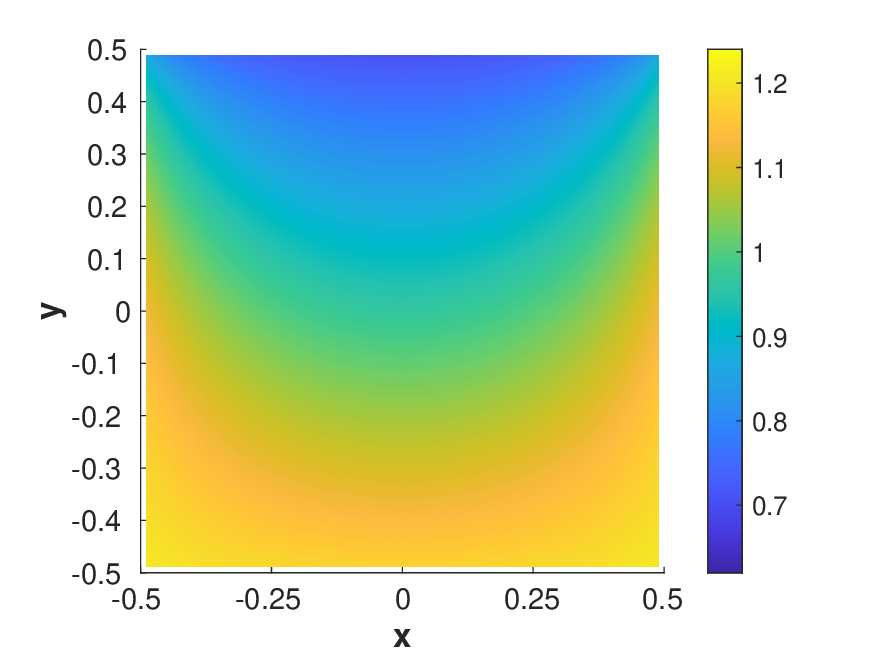}
    }
    \subfigure[$\epsilon = 10^{-2}$]{
    \includegraphics[width=0.31\textwidth, trim=20 2 50 20, clip]{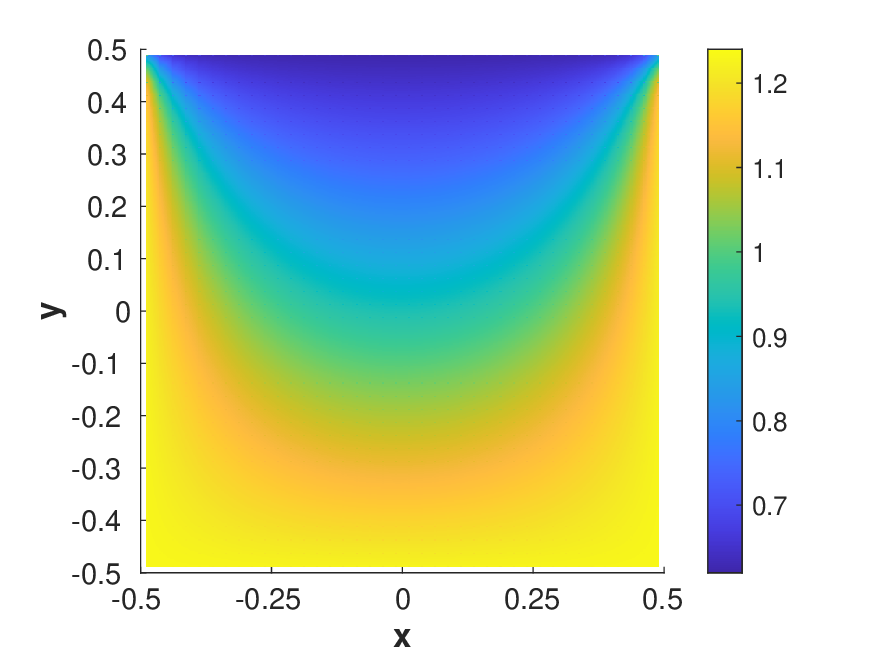}
    } 
    \subfigure[$\epsilon = 10^{-4}$]{
    \includegraphics[width=0.31\textwidth, trim=20 2 50 20, clip]{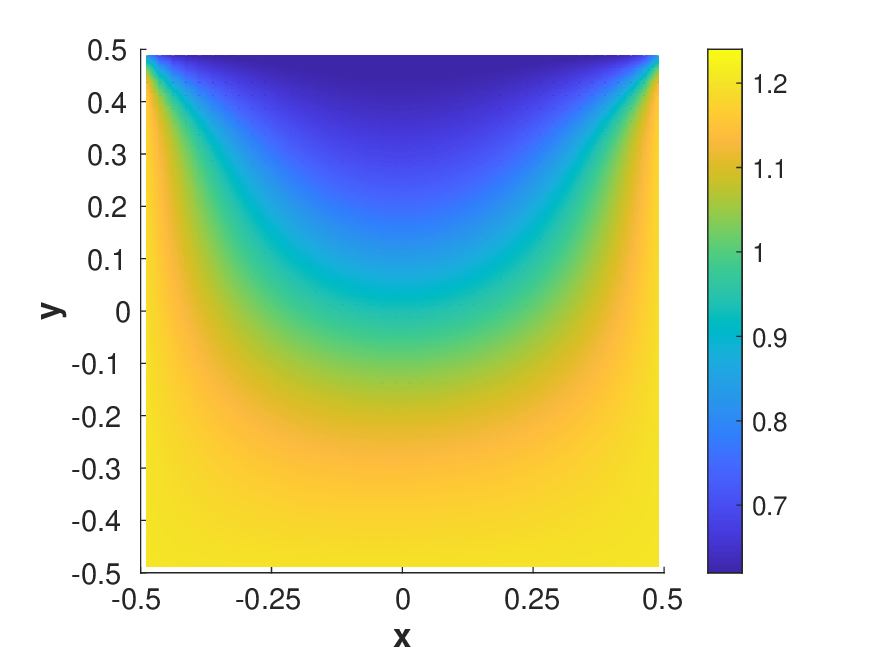}
    }
    \caption{The distributions of density for heat transfer in a cavity.}
\end{figure}

\begin{figure}[!ht] \label{fig:T_2D3D_BGK}
    \centering
    \subfigure[$\epsilon = 1$]{
    \includegraphics[width=0.31\textwidth, trim=20 2 50 18, clip]{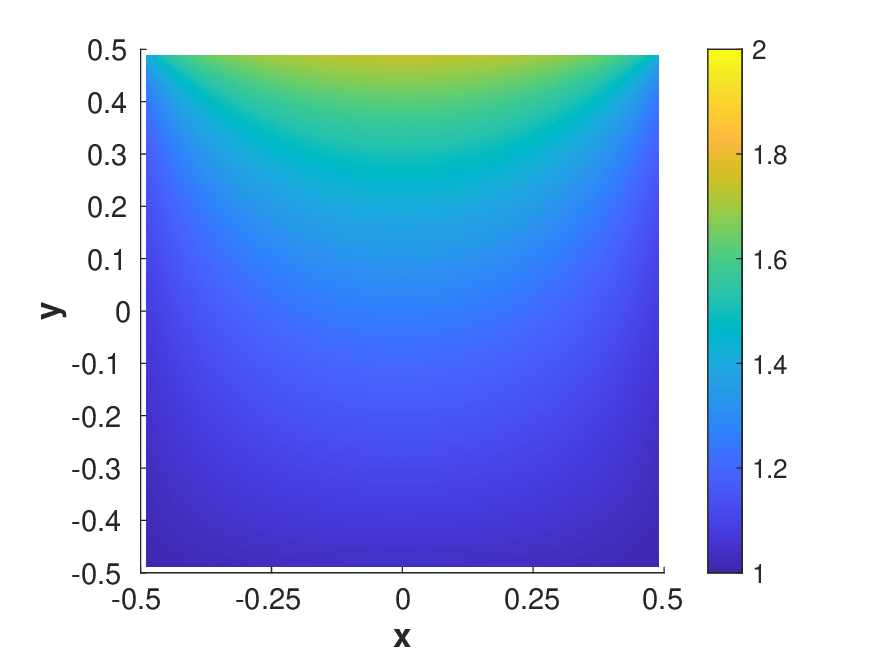}
    }
    \subfigure[$\epsilon = 10^{-2}$]{
    \includegraphics[width=0.31\textwidth, trim=20 2 50 18, clip]{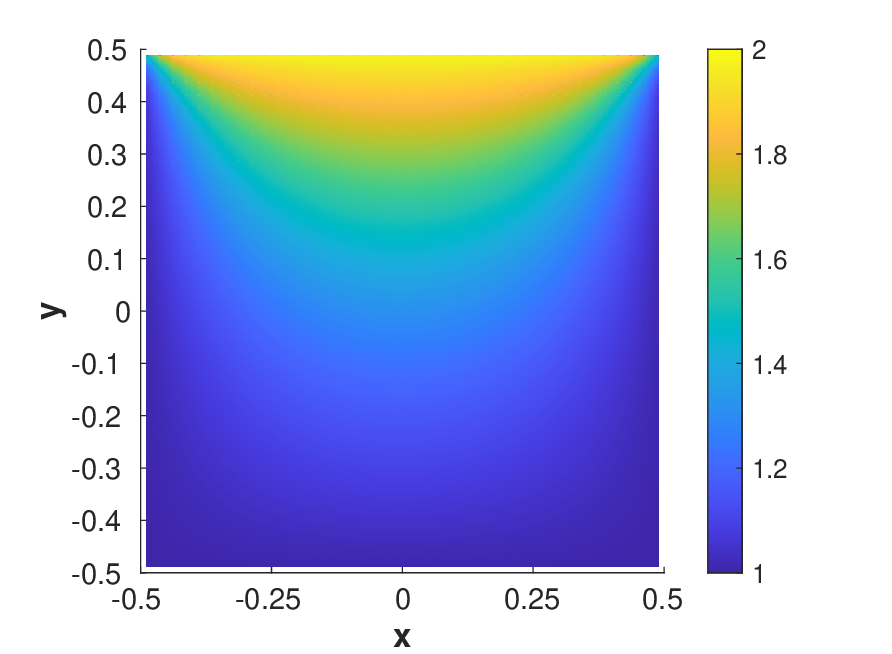}
    }
    \subfigure[$\epsilon = 10^{-4}$]{
    \includegraphics[width=0.31\textwidth, trim=20 2 50 18, clip]{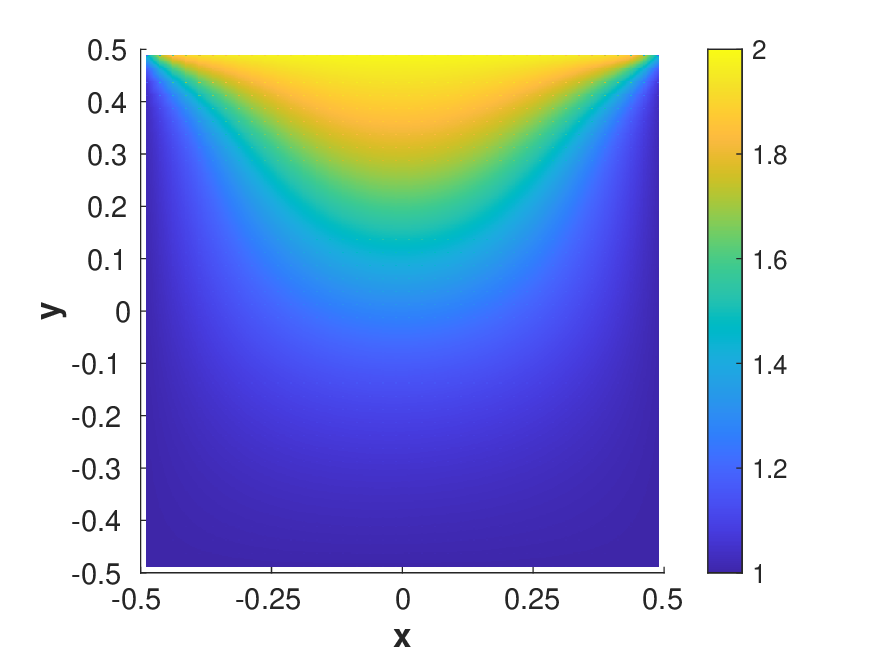}
    }
    \caption{The distributions of temperature for the heat transfer in a cavity.}
\end{figure}

Our second example for the BGK model is a 2D3V heat transfer problem where the domain is a unit square $\Omega = (-1/2,1/2) \times (-1/2,1/2)$, and the parameter $\nu$ is chosen as $4 \pi \rho$ in \eqref{eq:BGK_collis} to match the loss term of the binary collision operator. All the four boundaries of the domain are considered as solid wall with no velocities. The temperature of the top wall $(-1/2,1/2) \times \{1/2\}$ is set to be $2$, and all other walls have temperature $1$. A $40 \times 40$ uniform grid is adopted to discretize the spatial domain, and the three-dimensional velocity domain is truncated to $[-L,L]^3$ with $L=6$ which is discretized into a $20^3$ uniform grid. The density and temperature obtained by the second-order scheme are displayed in \Cref{fig:rho_2D3D_BGK} and \ref{fig:T_2D3D_BGK}, respectively. When $\epsilon$ becomes smaller, the jumps of density and temperature in the top-left and top-right corners are more obvious due to the discontinuities of temperature in boundary conditions.

\begin{figure}[!ht] \label{fig:rate_2D3V}
    \centering
    \subfigure[SGS-FP (1st-order)]{\label{fig:rate_2D3V_FP}
    \includegraphics[width=0.31\textwidth, trim=5 0 35 20, clip]{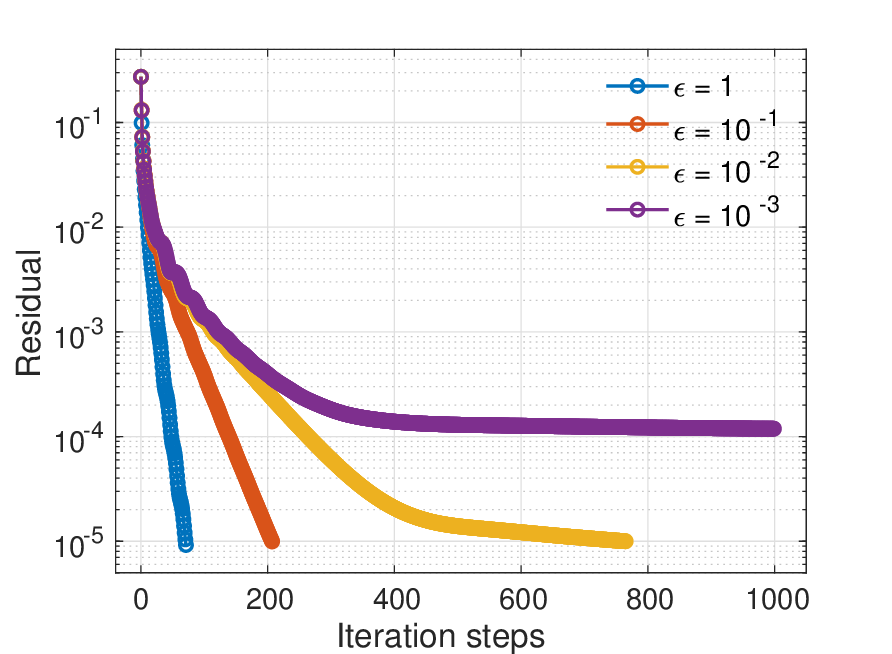}
    } 
    \subfigure[SGS-PFP (1st-order)]{\label{fig:rate_2D3V_MM_1st}
    \includegraphics[width=0.31\textwidth, trim=5 0 35 20, clip]{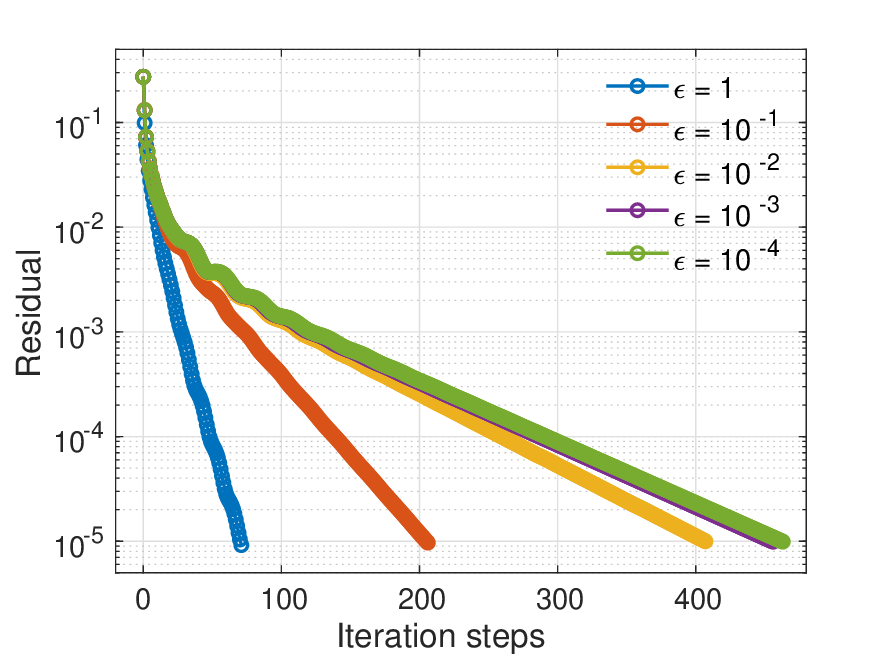}
    } 
    \subfigure[SGS-PFP (2nd-order)]{\label{fig:rate_2D3V_MM_2nd}
    \includegraphics[width=0.303\textwidth, trim=5 0 35 20, clip]{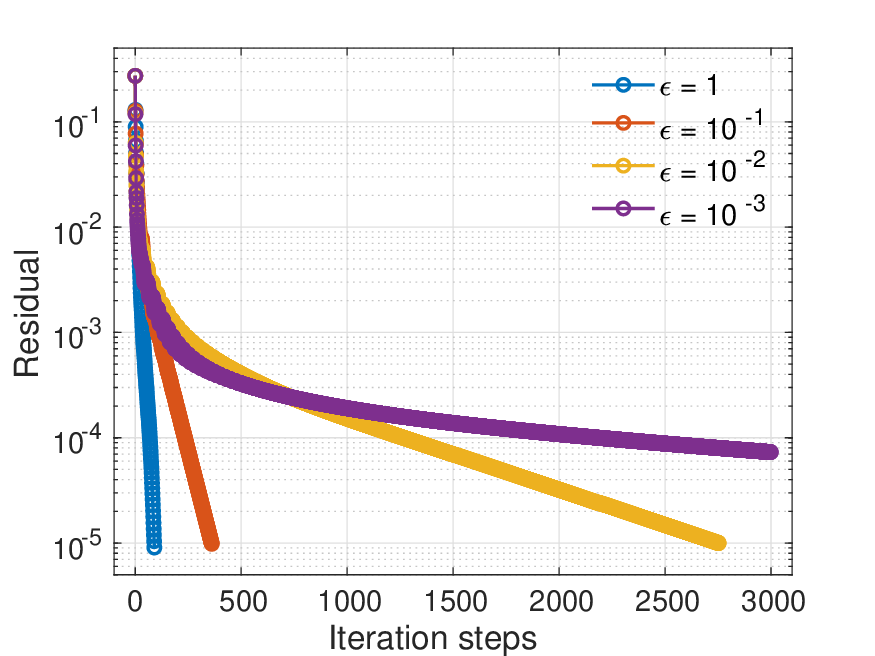}
    }
    \caption{The convergence of iterative methods for heat transfer in a cavity.}
\end{figure}

We again focus mainly on the efficiency of iterative methods. Similar to the one-dimensional case (see \Cref{fig:rate_1D1V}), it can be seen from \Cref{fig:rate_2D3V} that more iterations are needed as $\epsilon$ declines. Theoretically, if the nonlinear system \eqref{eq:inner_eqn} is solved exactly, the residual of SGS method \eqref{eq:SGS_BGK} does not depend on the algorithm to solve \eqref{eq:inner_eqn}, meaning that the Figures \ref{fig:rate_2D3V_FP} and \ref{fig:rate_2D3V_MM_1st} should be identical. But in practice, the preconditioned method can achieve a much better accuracy when solving \eqref{eq:inner_eqn}, which then leads to a much lower number of outer iterations when $\epsilon$ is small. In particular, when $\epsilon = 10^{-3}$, the residual of the SGS-FP method hardly decreases after $600$ iterations. For second-order methods, again we observe a decline in the efficiency of the SGS method, which needs to be improved by the multigrid technique.

\begin{figure}[!ht]  \label{fig:inner_1st_2D3V}
    \centering
    \subfigure[$\epsilon = 1$]{
    \includegraphics[width=0.4\textwidth, trim=2 2 35 15, clip]{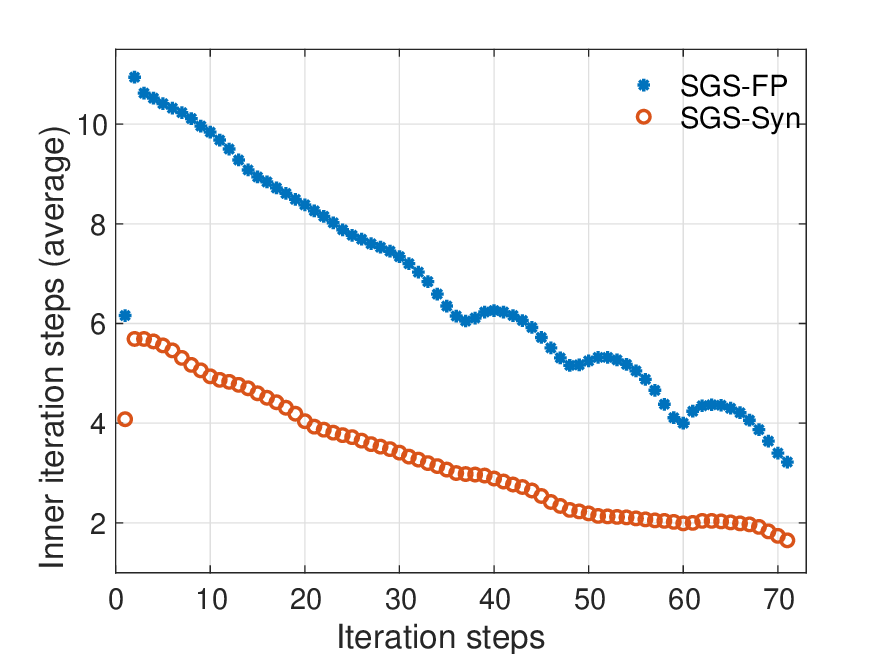}
    } \quad
    \subfigure[$\epsilon = 10^{-1}$]{
    \includegraphics[width=0.4\textwidth, trim=2 2 35 15, clip]{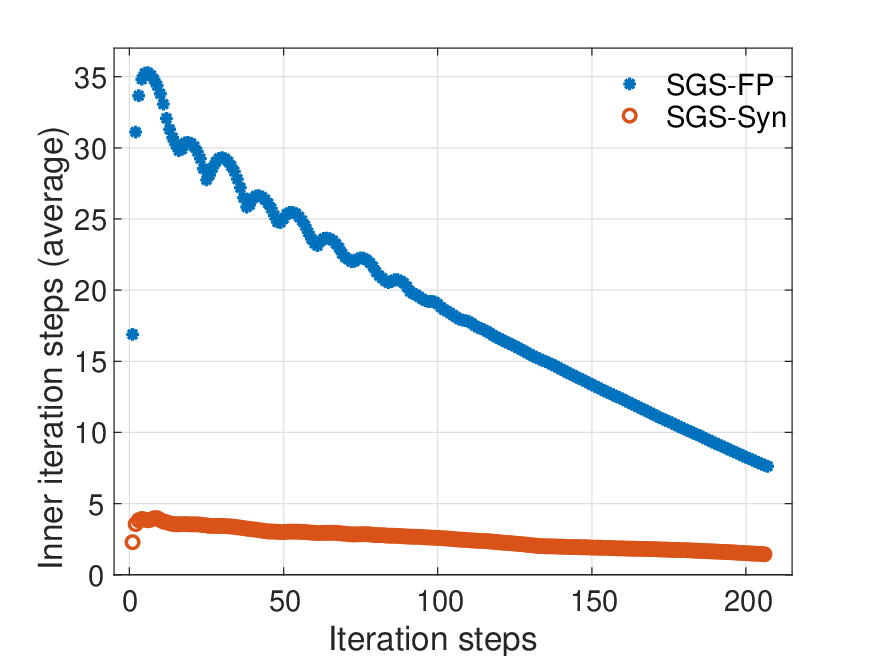}
    } \\
    \subfigure[$\epsilon = 10^{-2}$]{
    \includegraphics[width=0.4\textwidth, trim=2 2 35 15, clip]{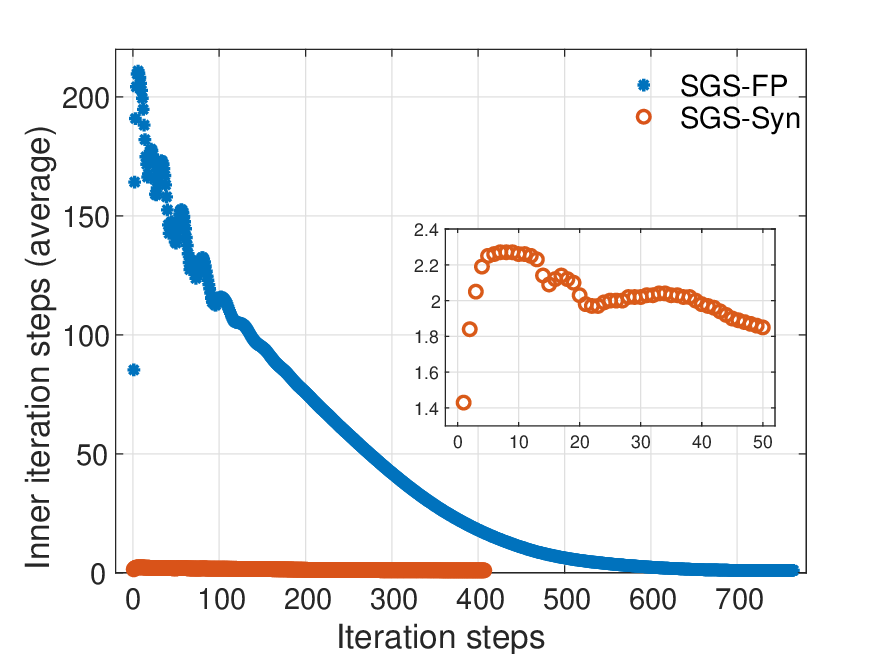}
    } \quad
    \subfigure[$\epsilon = 10^{-3}$]{\label{fig:inner_1st_2D3V_3}
    \includegraphics[width=0.4\textwidth, trim=2 2 32 15, clip]{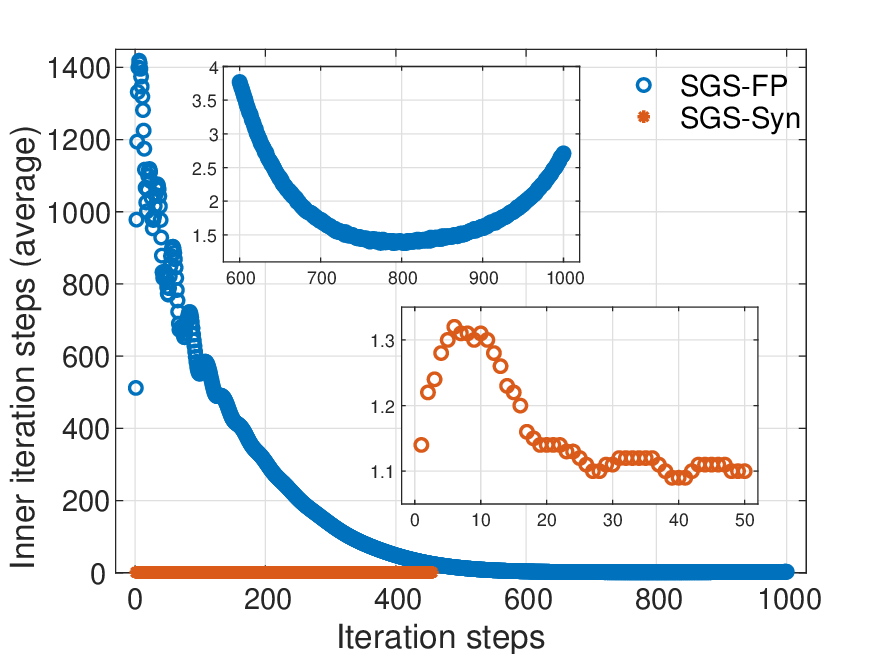}
    } \\
    \caption{The average inner iterations of iterative methods (1st-order) for heat transfer in a cavity.}
\end{figure}

The average numbers of inner iterations are plotted in \Cref{fig:inner_1st_2D3V}. The meanings of the axes are the same as those in \Cref{fig:inner_1st_1D1V}. Compared with the one-dimensional case, the difference between two methods becomes even sharper for the 3D velocity domain. While the fixed-point iteration requires more iterations for smaller $\epsilon$, the number of iterations for the preconditioned method decreases with $\epsilon$. For $\epsilon = 10^{-3}$, the inset in \Cref{fig:inner_1st_2D3V_3} shows that the curve  for SGS-FP starts to rise after 80 outer iterations, indicating that a lower threshold for the fixed-point iteration needs to be set to meet the termination condition of the SGS method.

\begin{figure}[!ht] \label{fig:mg_2D3V}
    \centering
    \subfigure[1st-order scheme]{
    \includegraphics[width=0.4\textwidth, trim=2 2 40 20, clip]{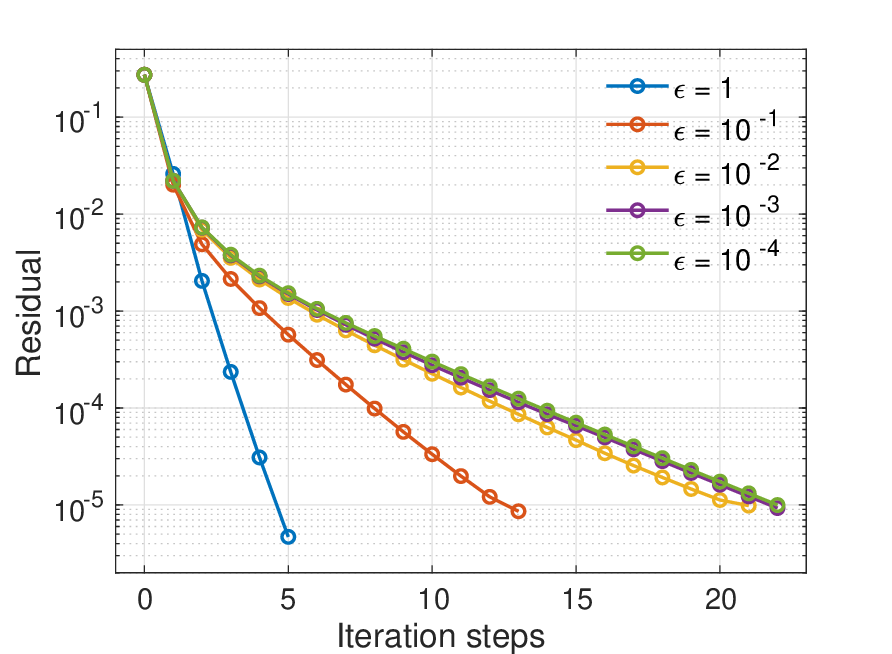}
    } \quad
    \subfigure[2nd-order scheme]{\label{fig:mg_2D3V_2nd}
    \includegraphics[width=0.41\textwidth, trim=2 2 28 20, clip]{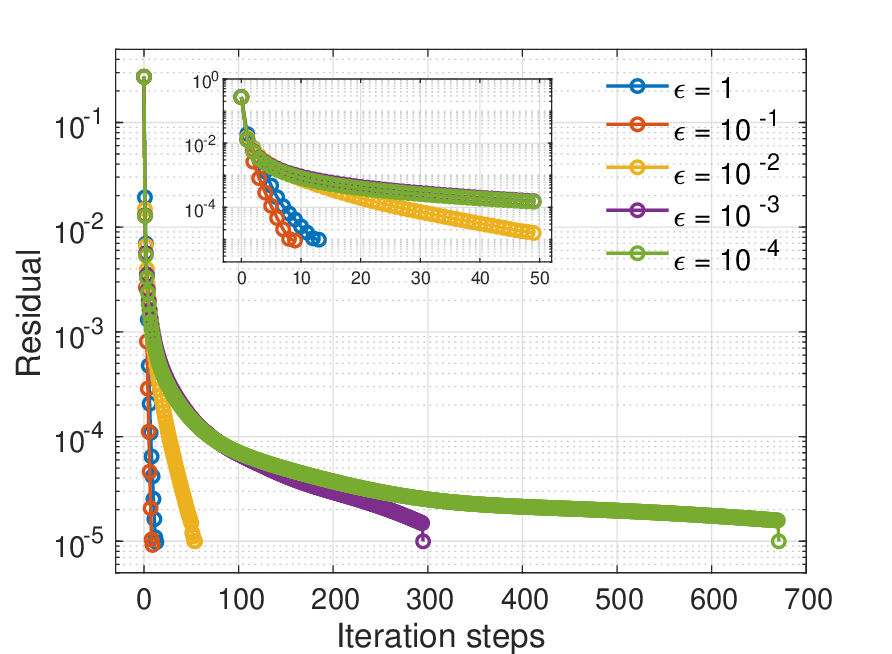}
    }
    \caption{The convergence of multigrid method for heat transfer in a cavity.}
\end{figure}

\Cref{fig:mg_2D3V} presents the convergence of the MG-SGS-PFP method. Again, the multigrid technique considerably reduces the number of iterations. A noteworthy phenomenon is that for the second-order scheme (\Cref{fig:mg_2D3V_2nd}) with $\epsilon = 10^{-3}$ and $10^{-4}$, the last iteration provides a sudden drop of the residual. It is observed that when this occurs, the V-cycle of the multigrid method does not coarsen the grid to the coarsest level. The residual at an intermediate grid is already reduced below the threshold by the presmoothing operator.

The computational time of SGS-FP and SGS-PFP is tabulated in \Cref{tab:time_2D3V} where ``$-$'' represents that the computation is terminated before the stopping criterion is met, since the time is longer than $48$ hours (the same for \Cref{tab:time_1D2V} and \Cref{tab:time_cavity_flow}). For this 2D3V problem, the SGS-PFP method outperforms the SGS-FP method in all cases, and the multigrid method still takes a significant role, especially for small values of $\epsilon$.
\begin{table}[!ht] \label{tab:time_2D3V}
    \caption{Computing time for heat transfer in a cavity. (measure time by the second).}
    \footnotesize
    \centering
    \begin{tabular}{cccccccc}
        \hline 
        \hline 
        Order & & & $\epsilon = 1$ & $\epsilon = 10^{-1}$ & $\epsilon = 10^{-2}$ & $\epsilon = 10^{-3}$ & $\epsilon = 10^{-4}$ \\
        \hline
        \multirow{3}{*}{First} & SGS-FP & & $2461.28$ & $18921.84$ & $162168.90$ & -- & -- \\
& SGS-PFP & & $1876.55$ & $4684.12$ & $6144.21$ & $5490.00$ & $5502.58$ \\
& MG-SGS-PFP & & $425.06$ & $1030.87$ & $1090.24$ & $971.69$ & $945.81$ \\
        \hline
        \multirow{2}{*}{Second} & SGS-PFP & & $2934.51$ & $9504.34$ & $40256.89$ & $171357.22$ & -- \\
& MG-SGS-PFP & & $2255.97$ & $1542.45$ & $6122.38$ & $30216.12$ & $67737.60$ \\
        \hline
        \hline 
    \end{tabular}
\end{table}

\subsection{Binary collision operator}

We now perform simulations for the binary collision operator $\mathcal{Q}[f,f]$ defined in \eqref{eq:collision}. For simplicity, we only consider the case of Maxwell molecules with $B(\cdot,\cdot)$ being a constant $B \equiv 1$, and in \eqref{eq:Pf}, we choose the parameter $\nu$ to be $\nu = 2\pi\rho$ for the models with two-dimensional velocity space. By such a choice, the loss term of the BGK operator matches the loss term of the binary collision operator. To discretize the collision term \eqref{eq:collision}, the Fourier spectral method \cite{Mouhot2006, Dimarco2014}, one of the most popular methods, is adopted. The implementation assumes that the distribution function has a compact support in $B(0,R)$, and the velocity domain is truncated to $[-L, L]^d$. In the following examples, we choose
\begin{equation*}
L = \frac{3\sqrt{2} + 1}{2} R, \qquad R = 3.
\end{equation*}
We have also applied the steady-state preserving method \eqref{eq:discrete_Q} to guarantee that the  Maxwellian $\mathcal{M}$ given in \eqref{eq:discrete_Gaussian} is in the kernel of the discrete collision term.

\subsubsection{Heat transfer between parallel plates}

\begin{figure}[!ht]  \label{fig:rho_T_1D2V}
    \centering
    \subfigure[Density]{
    \includegraphics[width=0.41\textwidth, trim=2 2 30 15, clip]{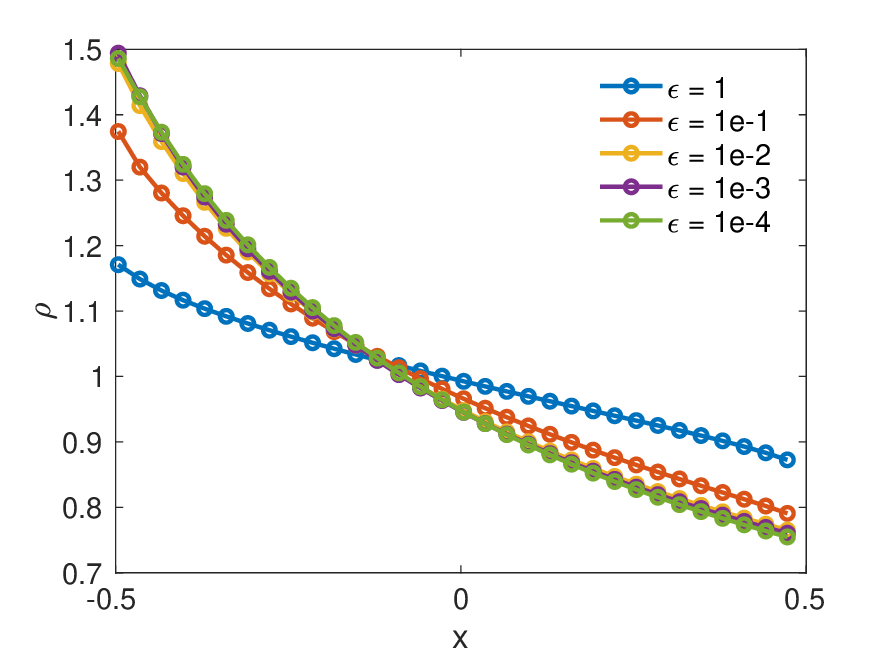}
    } \quad
    \subfigure[Temperature]{
    \includegraphics[width=0.4\textwidth, trim=2 2 30 15, clip]{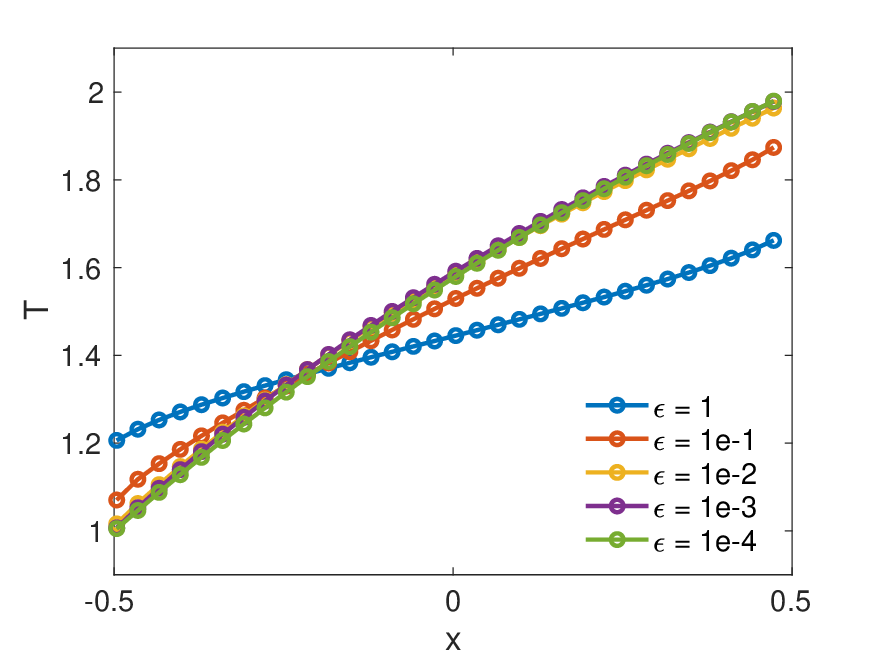}
    }
    \caption{The distributions of density and temperature for heat transfer between parallel plates.}
\end{figure}

Our first example is a 1D2V model of heat transfer between two parallel plates. The spatial domain is discretized on a uniform grid with $\Delta x = 1/128$, and the grid size of velocity domain is set to be $32 \times 32$. The density and temperature plots are given in \Cref{fig:rho_T_1D2V}, from which we again observe higher temperature jump and $\epsilon$ gets smaller. Note that when $\epsilon$ tends to zero, the temperature profile does not approach a straight line due to the dependence of the viscosity coefficient on $\rho$. This differs from the BGK case with a constant collision frequency as shown in \Cref{fig:rho_T_1D1V}.

\begin{figure}[!ht] 
    \centering 
    \subfigure[SGS-FP (1st-order)]{ \label{fig:rate_1D2V_FP}
    \includegraphics[width=0.3\textwidth, trim=2 2 35 20, clip]{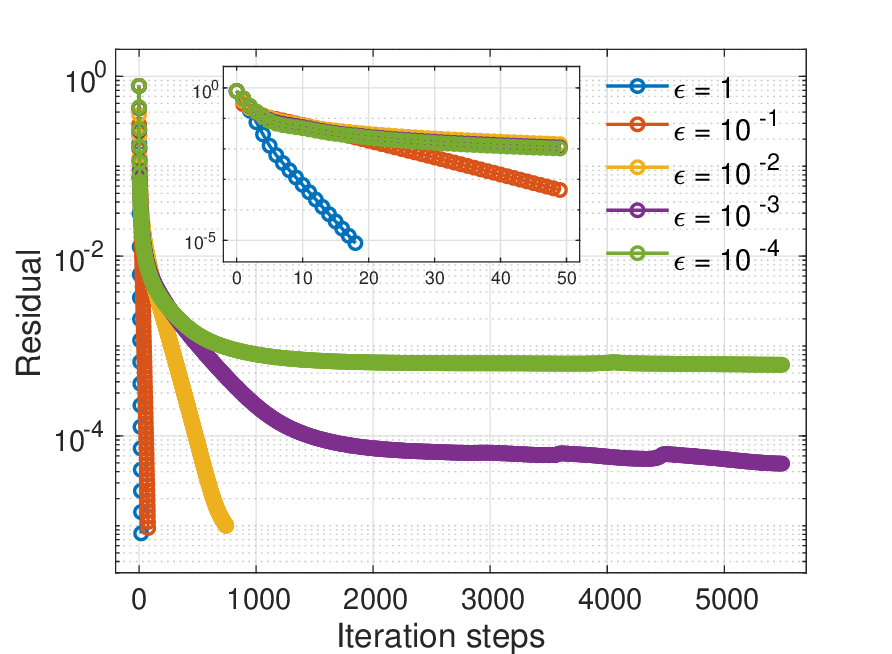}
    }
    \subfigure[SGS-PFP (1st-order)]{ \label{fig:rate_1D2V_MM_1}
    \includegraphics[width=0.3\textwidth, trim=2 2 35 20, clip]{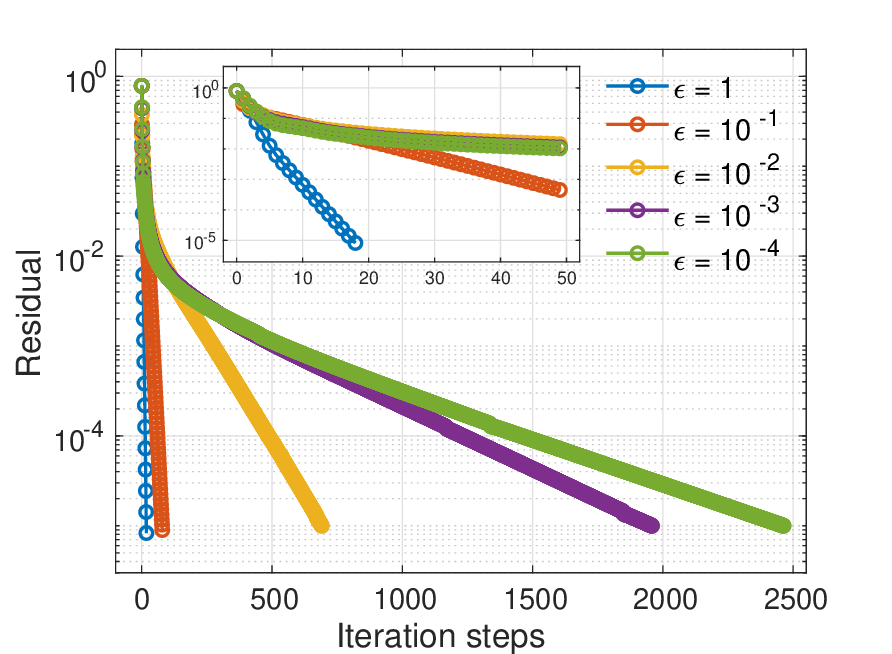}
    } 
    \subfigure[SGS-PFP (2nd-order)]{ \label{fig:rate_1D2V_MM_2}
    \includegraphics[width=0.3\textwidth, trim=2 2 35 20, clip]{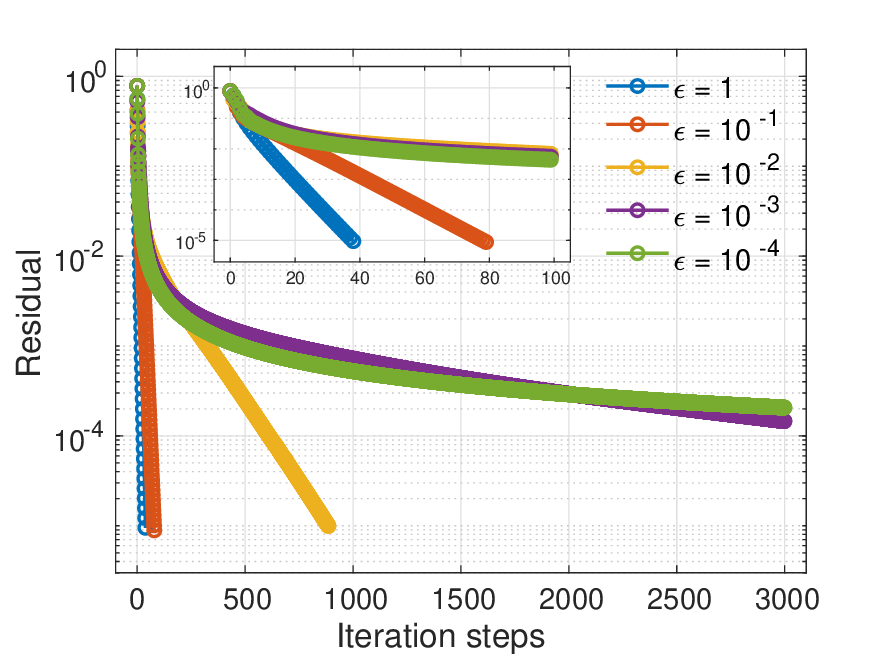}
    }
    \caption{The convergence of iterative methods for heat transfer between parallel plates.}
    \label{fig:rate_1D2V}
\end{figure}

\begin{figure}[!ht] \label{fig:inner_1D2V}
    \centering
    \subfigure[$\epsilon = 1$]{
    \includegraphics[width=0.4\textwidth, trim=2 2 35 15, clip]{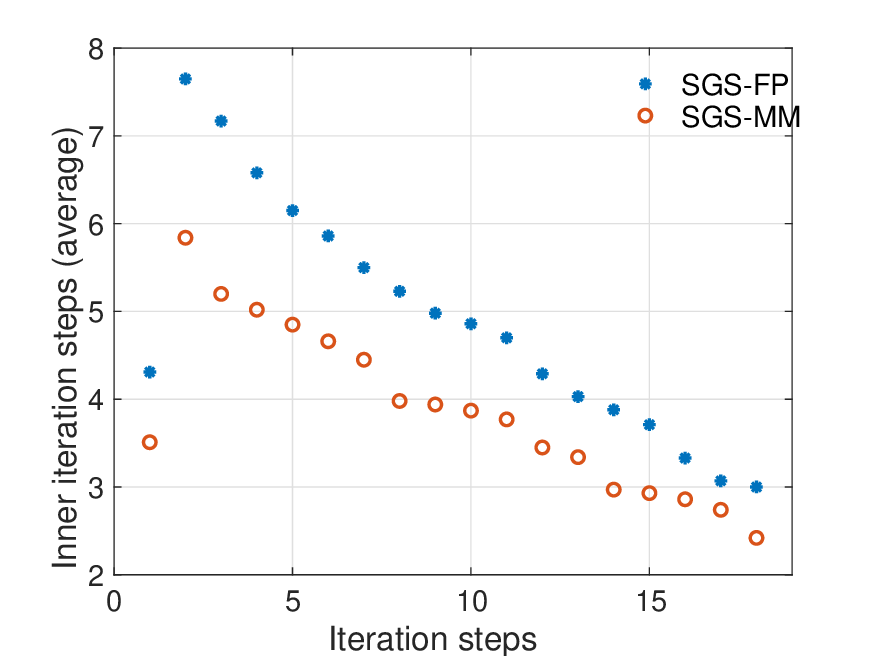}
    } \quad
    \subfigure[$\epsilon = 10^{-1}$]{
    \includegraphics[width=0.4\textwidth, trim=2 2 35 15, clip]{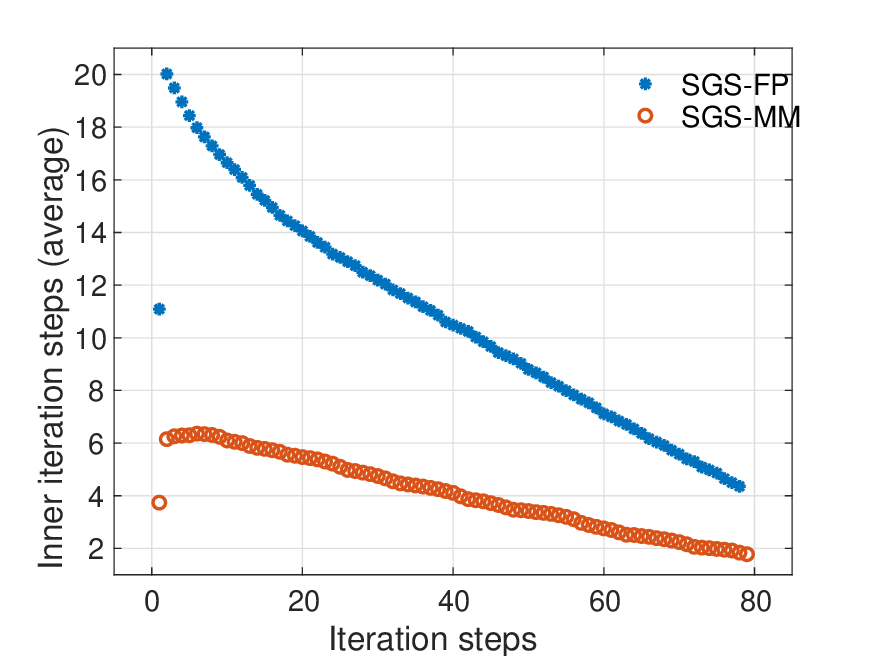}
    } \\
    \subfigure[$\epsilon = 10^{-2}$]{
    \includegraphics[width=0.4\textwidth, trim=2 2 35 15, clip]{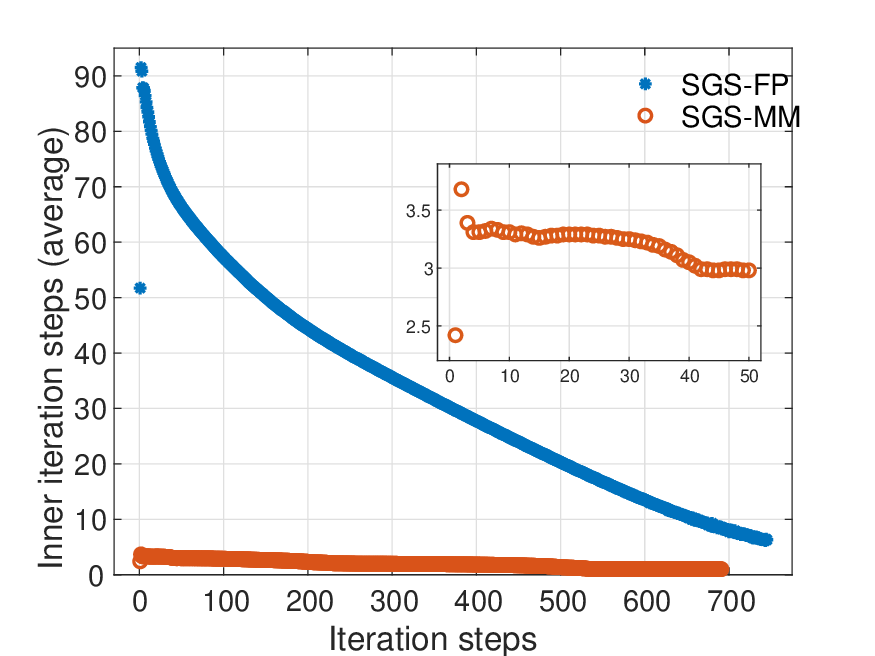}
    } \quad
    \subfigure[$\epsilon = 10^{-3}$]{
    \includegraphics[width=0.4\textwidth, trim=2 2 32 15, clip]{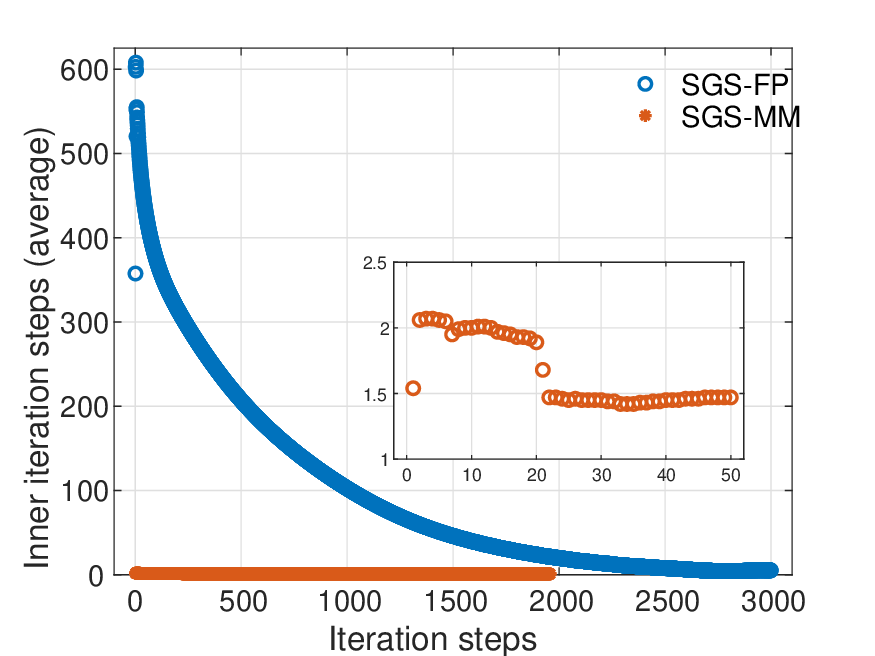}
    } \\
    \caption{The average inner iterations of iterative methods (1st-order) for heat transfer between parallel plates.}
\end{figure}

The evolution of residuals for SGS-based methods are plotted in \Cref{fig:rate_1D2V}. Note that the source iteration is not studied for this example since the computation of the binary collision operator is much more time-consuming than the BGK operator. The general behavior is similar to the BGK case. We again observe that the quality of inner iterations actually affects the number of outer iterations, and for smaller $\epsilon$, the system is closer to elliptic, requiring more iterations to reach the threshold. The average numbers of inner iterations presented in \Cref{fig:inner_1D2V} further illustrate the efficiency of SGS-PFP when $\epsilon$ is small, whose average numbers of inner iterations do not increase as $\epsilon$ gets small.


\begin{figure}[!ht] \label{fig:mg_1D2V}
    \centering
    \subfigure[1st-order scheme]{
    \includegraphics[width=0.4\textwidth, trim=2 2 40 20, clip]{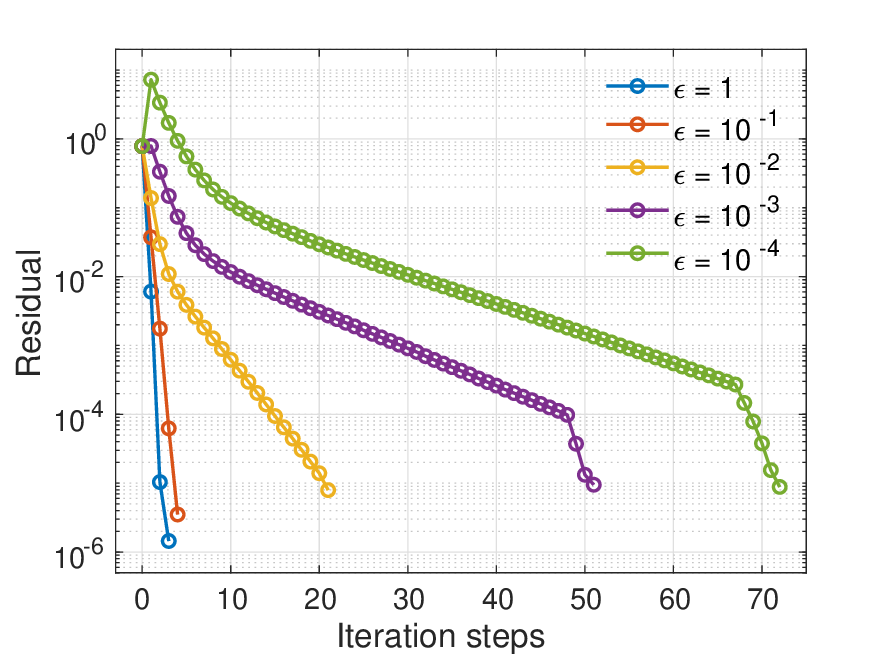}
    } \quad
    \subfigure[2nd-order scheme]{
    \includegraphics[width=0.4\textwidth, trim=2 2 40 20, clip]{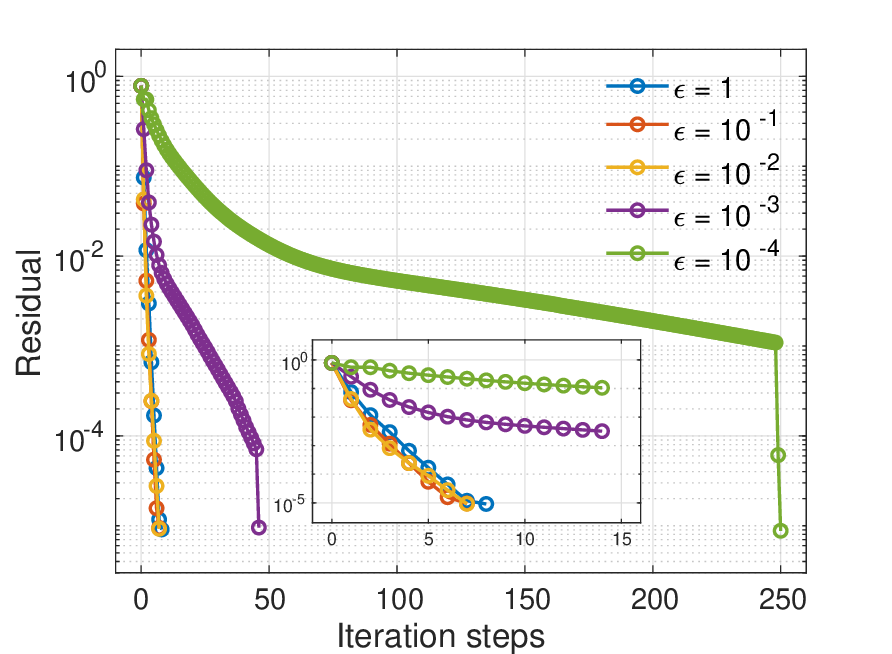}
    }
    \caption{The convergence of multigrid method for heat transfer between parallel plates.}
\end{figure}

The multigrid technique again takes a crucial role in enhancing the performance, as demonstrated in \Cref{fig:mg_1D2V}. For $\epsilon = 10^{-3}$ and $10^{-4}$, rapid declines of the residual occur again in the last few steps, and in these steps, the multigrid method does not reach the coarsest level. \Cref{tab:time_1D2V} lists the computational time of the above iterative methods. The results confirm that both the preconditioning and the multigrid technique can significantly reduce the computational time, despite the overhead in maintaining the data structure and copying the data between coarse and fine grids in the implementation of multigrid.
\begin{table}[!ht] \label{tab:time_1D2V}
    \caption{Computing time for heat transfer between parallel plates (measure time by the second).}
    \footnotesize
    \centering
    \begin{tabular}{cccccccc}
        \hline 
        \hline 
        Order & & & $\epsilon = 1$ & $\epsilon = 10^{-1}$ & $\epsilon = 10^{-2}$ & $\epsilon = 10^{-3}$ & $\epsilon = 10^{-4}$ \\
        \hline
        \multirow{3}{*}{First} & SGS-FP & & $18.39$ & $105.41$ & $1855.94$ & $73320.70$ & -- \\
& SGS-PFP & & $17.05$ & $97.34$ & $646.96$ & $1582.93$ & $1968.14$ \\
& MG-SGS-PFP & & $12.73$ & $25.59$ & $112.09$ & $252.94$ & $354.35$ \\
        \hline
        \multirow{2}{*}{Second} & SGS-PFP & & $40.99$ & $97.15$ & $831.70$ & $5572.82$ & $28495.87$ \\
& MG-SGS-PFP & & $69.82$ & $55.77$ & $78.17$ & $466.94$ & $2512.15$ \\
        \hline
        \hline 
    \end{tabular}
\end{table}

\subsubsection{Lid-driven cavity flow}
\begin{figure}[!ht] \label{fig:rho_flow}
    \centering
    \subfigure[$\epsilon = 1$]{
    \includegraphics[width=0.307\textwidth, trim=10 2 48 20, clip]{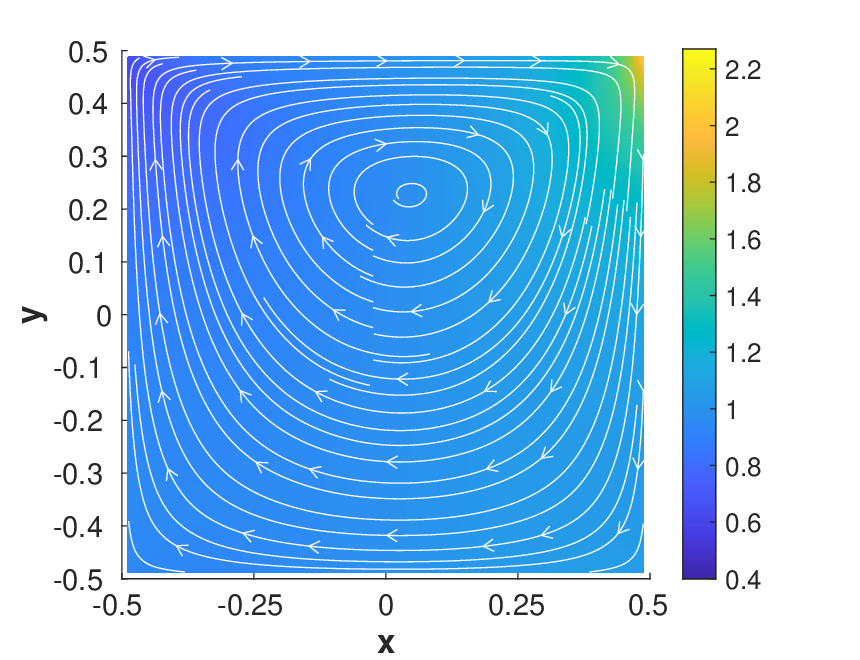}
    }
    \subfigure[$\epsilon = 10^{-2}$]{
    \includegraphics[width=0.31\textwidth, trim=20 2 50 20, clip]{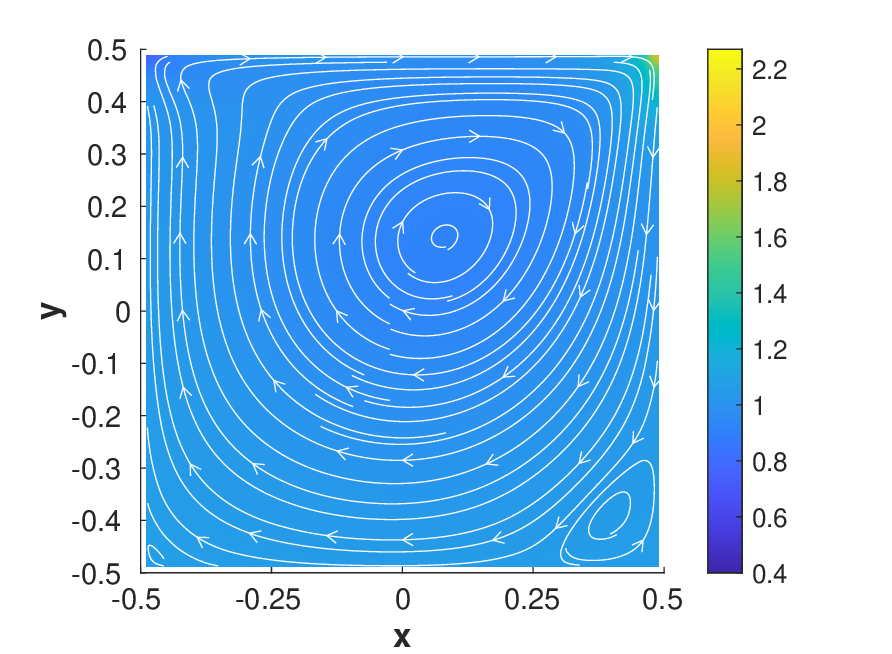}
    }
    \subfigure[$\epsilon = 10^{-4}$]{
    \includegraphics[width=0.31\textwidth, trim=20 2 50 20, clip]{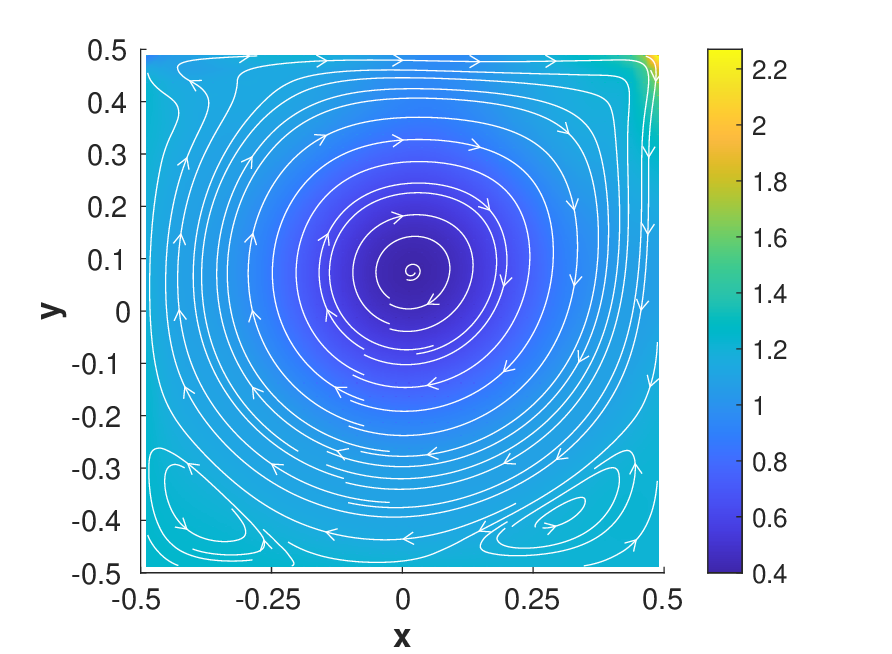}
    }    
    \caption{The distributions of density and the streamlines of macroscopic velocity for lid-driven cavity flow.}
\end{figure}

Lid-driven cavity flow is a classical benchmark problem in fluid dynamics, offering insights into fluid behavior under controlled conditions and serving as a basis for understanding more complex flow phenomena. Here we consider a 2D2V model where gas is confined in a unit square cavity $\Omega = (-1/2,1/2) \times (-1/2,1/2)$ with a lid at top boundary $y=1/2$ that moves at a constant velocity $\bs{U}_w = (1,0)^T$ while the other boundaries remain stationary. All the walls have the same temperature $1$. The spatial domain and velocity domain are discretized uniformly with $40 \times 40$ and $32 \times 32$ grids, respectively. \Cref{fig:rho_flow} displays the numerical results of density and the streamlines of macroscopic velocity. As $\epsilon$ gets smaller, it can be observed that the position of the main vortex center gradually shifts from the upper right to the center, and when $\epsilon = 10^{-2}$, a secondary vortex appears in the bottom right corner whose rotation is opposite to the main vortex. Furthermore, for $\epsilon = 10^{-4}$, the coverage area of the vortex in the bottom right corner gradually expands and the center moves to the middle, and there is another secondary vortex in the bottom left corner. The above results agree with the classical fluid mechanics, where small Knudsen numbers correspond to high Reynolds numbers.

\begin{figure}[!ht] \label{fig:T_flow}
    \centering
    \subfigure[$\epsilon = 1$]{
    \includegraphics[width=0.31\textwidth, trim=20 2 50 20]{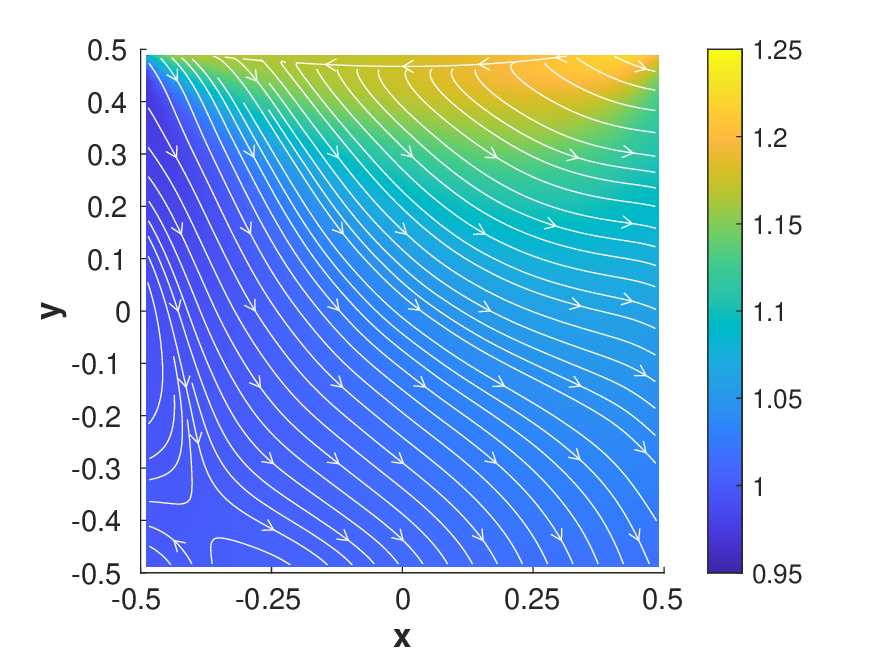}
    }
    \subfigure[$\epsilon = 10^{-2}$]{
    \includegraphics[width=0.31\textwidth, trim=20 2 50 20]{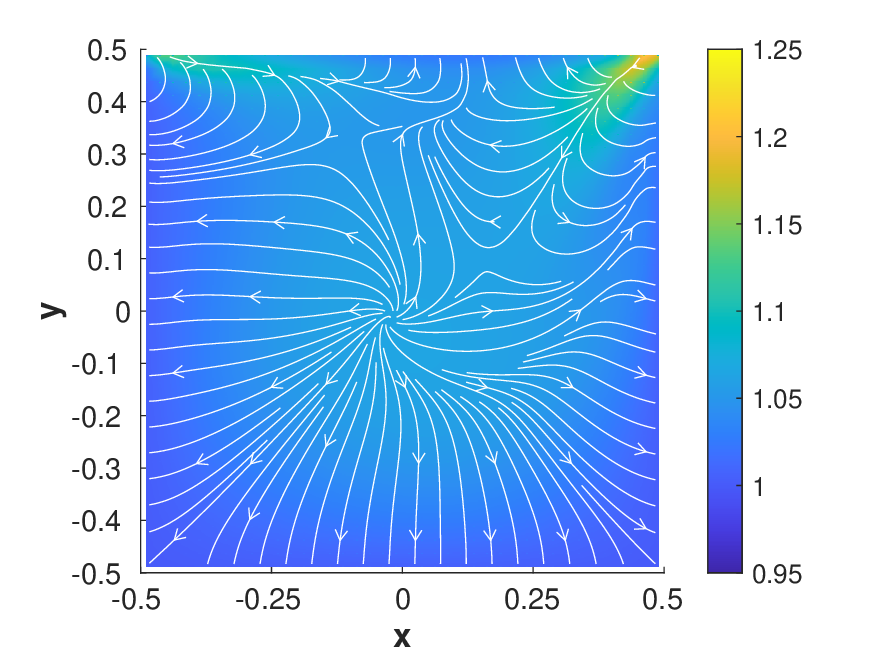}
    }
    \subfigure[$\epsilon = 10^{-4}$]{
    \includegraphics[width=0.31\textwidth, trim=20 2 50 20]{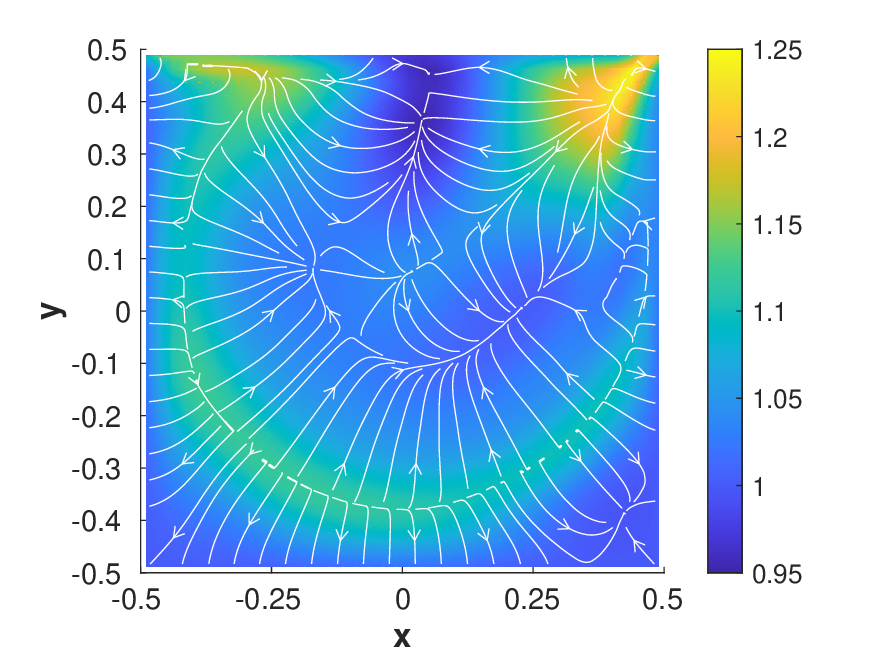}
    }    
    \caption{The distributions of temperature and the streamlines of heat flux for lid-driven cavity flow.}
\end{figure}

In \Cref{fig:T_flow}, the distributions of temperature and the streamlines of heat flux are plotted where the heat flux is defined by 
\begin{equation*}
\bs{q}(\bs{x}) = \frac{1}{2} \int_{\mathbb{R}^2} (\bs{v} - \bs{U}(\bs{x})) |\bs{v} - \bs{U}(\bs{x})|^2 f(\bs{x}, \bs{v}) \,\mathrm{d}\bs{v}, \quad i = 1,2  
\end{equation*} 
where $\bs{U}$ is the macroscopic velocity defined in \eqref{eq:moment}. When $\epsilon = 1$, it is clear that the heat flows from hot area to cold area, which reflects a typical rarefaction effect. For $\epsilon = 10^{-2}$ and $\epsilon = 10^{-4}$, the Fourier flow generally holds in the central region of the computational domain, but in the area close to the top lid, anti-Fourier heat flows can still be observed.

\begin{figure}[!ht] \label{fig:rate_flow_SGS}
    \centering 
    \subfigure[SGS-FP (1st-order)]{
    \includegraphics[width=0.3\textwidth, trim=2 2 35 20, clip]{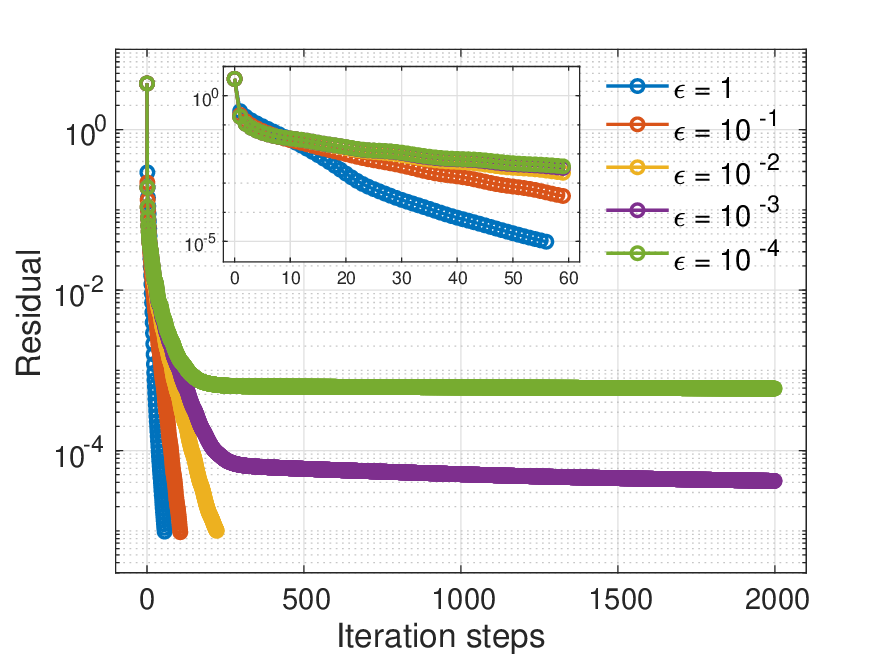}
    }
    \subfigure[SGS-PFP (1st-order)]{
    \includegraphics[width=0.3\textwidth, trim=2 2 35 20, clip]{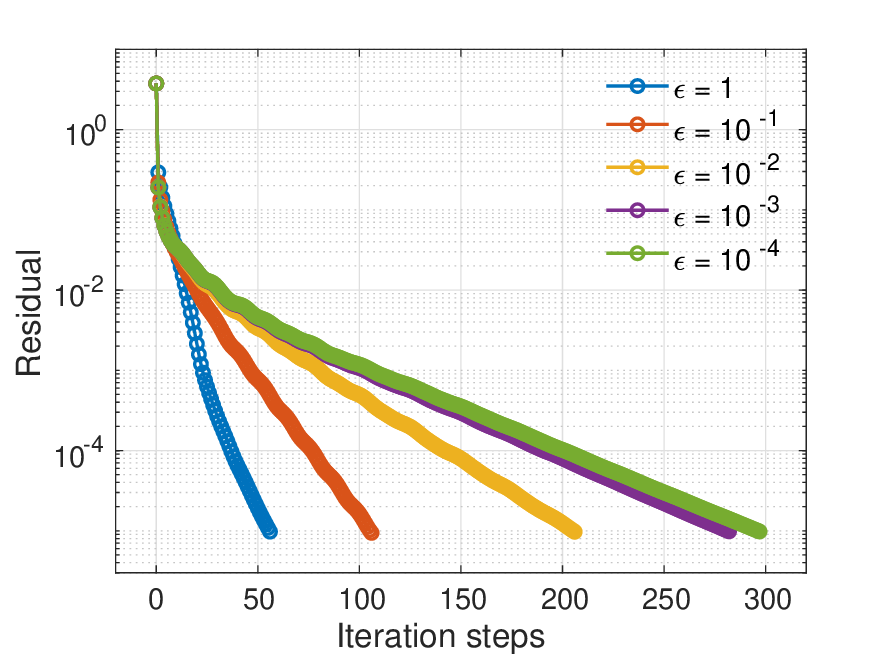}
    } 
    \subfigure[SGS-PFP (2nd-order)]{\label{fig:rate_flow_MM_2}
    \includegraphics[width=0.3\textwidth, trim=2 2 35 20, clip]{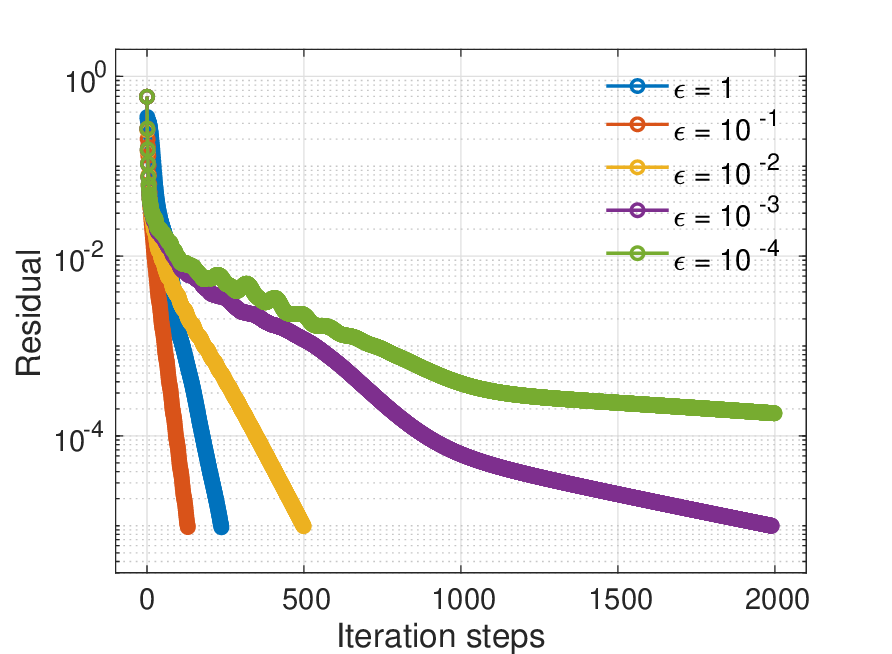}
    }
    \caption{The convergence of iterative methods for lid-driven cavity flow.}
\end{figure}

\begin{figure}[!ht] \label{fig:rate_flow_MG}
    \centering
    \subfigure[1st-order scheme]{
    \includegraphics[width=0.408\textwidth, trim=2 2 33 20, clip]{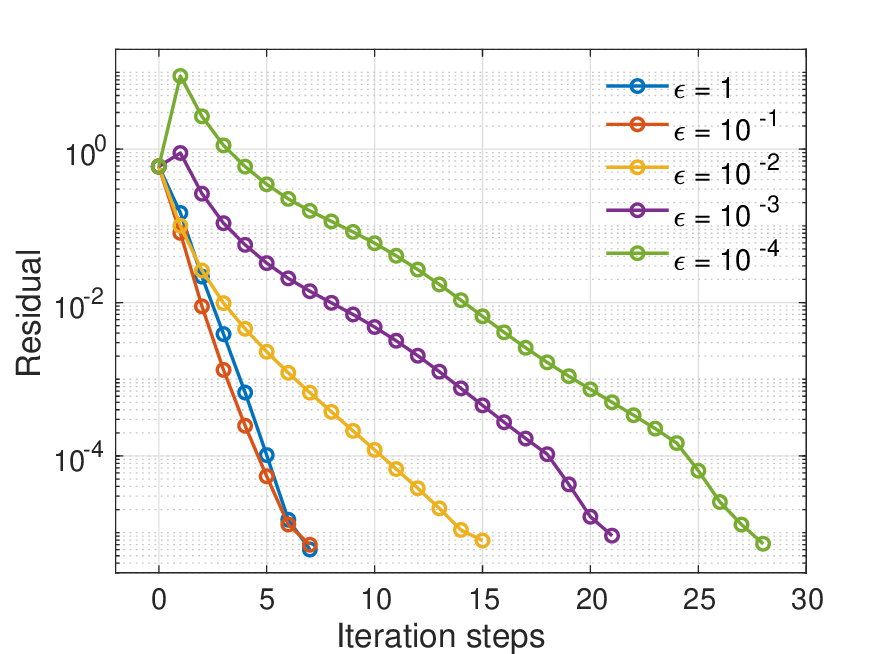}
    } \quad
    \subfigure[2nd-order scheme]{\label{fig:rate_flow_MG_2}
    \includegraphics[width=0.4\textwidth, trim=2 2 38 20, clip]{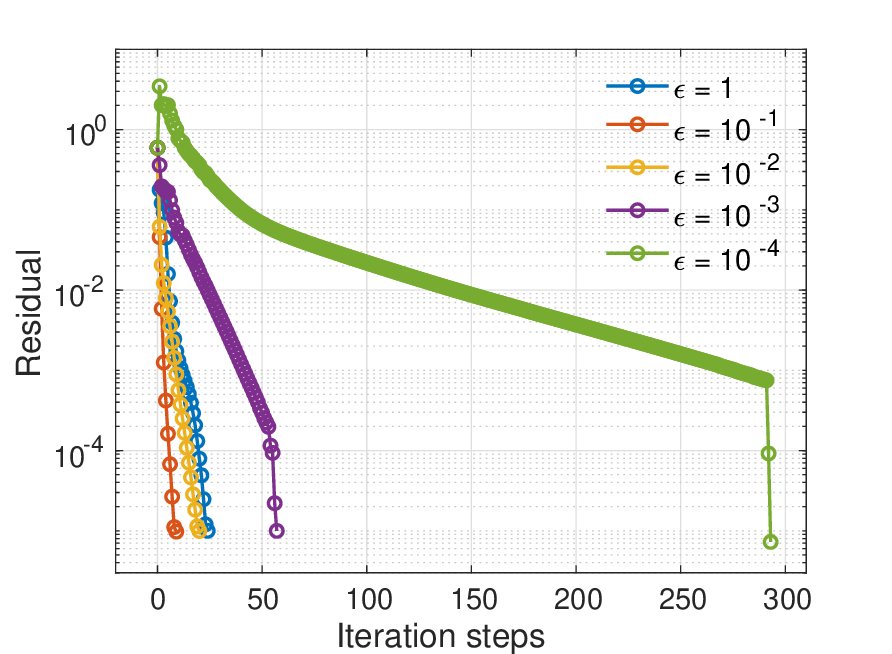}
    }
    \caption{The convergence of the MG-SGS-PFP method for lid-driven cavity flow.}
\end{figure}

\Cref{fig:rate_flow_SGS} and \Cref{fig:rate_flow_MG} show the decay of residual during the iteration, which resemble the plots in previous examples. One difference in this example is that in the second-order SGS-PFP method, the convergence rate for $\epsilon = 1$ is slower than that of $\epsilon = 10^{-1}$, which contradicts previous cases. One possible reason is that the structure of the solution changes significantly as $\epsilon$ changes from $1$ to $10^{-2}$, so that a direct comparison of convergence rates is less meaningful.

Finally, we present the computational time in \Cref{tab:time_cavity_flow}, which is consistent with the convegence results. SGS-PFP still runs fast than SGS-FP when $\epsilon \leqslant 10^{-2}$, and MG-SGS-PFP runs faster than the other two methods in all cases.
\begin{table}[!ht]  \label{tab:time_cavity_flow}
    \caption{Computing time for lid-driven cavity flow. (measure time by the second).}
    \footnotesize
    \centering
    \begin{tabular}{cccccccc}
        \hline 
        \hline 
        Order & & & $\epsilon = 1$ & $\epsilon = 10^{-1}$ & $\epsilon = 10^{-2}$ & $\epsilon = 10^{-3}$ & $\epsilon = 10^{-4}$ \\
        \hline
        \multirow{3}{*}{First} & SGS-FP & & $466.71$ & $1096.63$ & $5500.71$ & $136181.67$ & -- \\
& SGS-PFP & & $537.60$ & $1206.91$ & $2117.99$ & $2720.22$ & $2852.99$ \\
& MG-SGS-PFP & & $223.88$ & $268.42$ & $569.22$ & $794.72$ & $1045.59$ \\
        \hline
        \multirow{2}{*}{Second} & SGS-PFP & & $2203.16$ & $1296.05$ & $5131.12$ & $19151.96$ & $77945.26$ \\
& MG-SGS-PFP & & $1588.20$ & $625.80$ & $1552.71$ & $4478.99$ & $22127.55$ \\
        \hline
        \hline 
    \end{tabular}
\end{table}

\section{Conclusion} 
\label{sec:conclusion}
We developed a numerical solver for the steady-state Boltzmann equation discretized by the finite volume method for the spatial variable. The iterative method to solve the nonlinear system is based on the symmetric Gauss-Seidel method, and we have also coupled it with the multigrid method to accelerate the convergence. The major contribution of this paper is to design an efficient iterative method to solve the local problem on each grid cell by employing an efficient preconditioner. The computational time is shortened significantly by this technique when $\epsilon$ is small.

Currently, one remaining issue in the numerical solver is the effectiveness of the multigrid method. Despite a significant reduction of the computational cost by the multigrid technique, the number of iterations still grows as the grid is refined, causing a large number of iterations for small $\epsilon$. In \cite{Hu2016}, it is mentioned that different cycles may result in different performances when the multigrid method is applied to kinetic equations. This will be investigated more deeply in our future work. For second-order schemes, better prolongation and restriction operators are needed to improve the convergence rates, which is also included in our ongoing works.

\bibliographystyle{siamplain}
\bibliography{Outline.bib}

\end{document}